\documentclass{amsart}

\newcommand{\baton}[1]{\mathbb #1}
\newcommand{\C}{{\baton C}}
\newcommand{\E}{{\baton E}}

\newcommand{\R}{{\baton R}}
\newcommand{\T}{{\baton T}}
\newcommand{\Z}{{\baton Z}}

\newcommand{\CA}{{\mathcal A}}
\newcommand{\CB}{{\mathcal B}}
\newcommand{\CC}{{\mathcal C}}
\newcommand{\CD}{{\mathcal D}}

\newcommand{\CI}{{\mathcal I}}

\newcommand{\CP}{{\mathcal P}}

\newcommand{\CX}{{\mathcal X}}

\newcommand{\norm}[1]{\lVert #1\rVert}

\newcommand{\one}{{\boldsymbol 1}}
\newcommand{\inv}{^{-1}}
\newcommand{\wh}{\widehat}

\newtheorem{theorem}{Theorem}[section]

\newtheorem{proposition}[theorem]{Proposition}
\newtheorem{corollary}[theorem]{Corollary}

\newtheorem*{theorem*}{Theorem}
\newtheorem*{lemma*}{Lemma}
\newtheorem*{proposition*}{Proposition}
\newtheorem*{corollary*}{Corollary}
\newtheorem*{vdC}{van der Corput's Lemma}

\theoremstyle{definition}
\newtheorem{definition}[theorem]{Definition}
\newtheorem*{notation}{Notation}

\theoremstyle{remark}
\newtheorem{remark}{Remark}[section]
\theoremstyle{plain}

\newcommand{\bzero}{{\boldsymbol 0}}

\newcommand{\ua}{{\bf a}}
\newcommand{\ub}{{\bf b}}
\newcommand{\uc}{{\bf c}}

\newcommand{\uI}{{\bf I}}
\newcommand{\uJ}{{\bf J}}

\newcommand{\ulg}{{\underline g}}
\newcommand{\ulh}{{\underline h}}
\newcommand{\ulq}{{\underline q}}
\newcommand{\ulz}{{\underline z}}
\newcommand{\ulx}{{\underline x}}
\newcommand{\uly}{{\underline y}}

\newcommand{\ulgamma}{{\underline \gamma}}

\newcommand{\nnorm}[1]{\lvert\!|\!| #1|\!|\!\rvert}
\newcommand{\type}[1]{^{[#1]}}

\DeclareMathOperator{\HP}{HP}

\newcommand{\LimAv}[1]{\text{\upshape lim averages }\bigl(#1\bigr)}
\newcommand{\LimAvI}[2]{\text{\upshape lim 
averages}_{#1}\bigl(#2\bigr)}

\newcommand{\LimsupAv}[1]{\text{\upshape limsup}\left|\text{\upshape 
averages}(#1)\right|}
\newcommand{\LimsupAvI}[2]{\text{\upshape limsup}\left|\text{\upshape 
averages}_{#1}(#2)\right|}

\newcounter{truc}
\newcommand{\trucenumi}{\setcounter{truc}{\value{enumi}}}
\newcommand{\enumitruc}{\setcounter{enumi}{\value{truc}}}

\begin{document}
\title{Uniformity seminorms on $\ell^\infty$ and applications}
\author{Bernard Host and Bryna Kra}

\address{ Universit\'e Paris-Est, Laboratoire d'analyse et de
math\'ematiques appliqu\'ees, UMR CNRS 8050, 5 bd Descartes, 77454 Marne la
Vall\'ee Cedex 2, France}
\email{bernard.host@univ-mlv.fr}
\address{Department of Mathematics, Northwestern University,
2033 Sheridan Road,  Evanston, IL 60208-2730, USA}
\email{kra@math.northwestern.edu}
\thanks{The second author  was partially supported by
NSF grant DMS-0555250.}

\begin{abstract}
A key tool in recent advances in understanding arithmetic 
progressions and other patterns in subsets of the integers 
is certain norms or seminorms.  One example is the 
norms on $\Z/N\Z$ introduced by Gowers in his proof of 
Szemer\'edi's Theorem, used to detect uniformity 
of subsets of the integers.  Another example is the seminorms 
on bounded functions in a measure preserving system (associated to
the averages in Furstenberg's proof of Szemer\'edi's Theorem) 
defined by the authors.
For each integer $k\geq 1$, we 
define seminorms on $\ell^\infty(\Z)$ analogous
to these norms and seminorms.  
We study the correlation of these norms with certain 
algebraically defined sequences, which arise from evaluating a continuous 
function on the homogeneous space of a nilpotent Lie group 
on a orbit (the nilsequences).
Using these seminorms, we define a dual norm that acts as an 
upper bound for the correlation of a bounded sequence with a nilsequence. 
We also prove an 
 inverse theorem for the seminorms, 
showing how a bounded sequence correlates with a nilsequence.
 As applications, we derive several 
ergodic theoretic results, including a nilsequence 
version of the Wiener-Wintner ergodic theorem, a
nil version of a corollary to the spectral theorem, 
and a weighted multiple ergodic convergence theorem. 
\end{abstract}

\maketitle
\section{Introduction}

\subsection{Norms and seminorms}

In his proof of Szemer\'edi's Theorem, Gowers~\cite{Gowers} introduced 
norms for functions defined on $\Z/N\Z$ that count 
parallelepiped configurations and can be used to detect certain 
patterns (such as arithmetic progressions) 
in subsets of the integers. 
In~\cite{HK}, we defined seminorms on bounded measurable 
functions on a measure preserving system, 
that can be viewed as averages over parallelepipeds and 
use them to control the norm of multiple ergodic averages (such 
as one evaluated along arithmetic progressions).  
Although the original definitions were quite different, 
it turns out that the  Gowers norms and the ergodic seminorms 
 are almost the same object, but are defined on different spaces:
one on the space of functions on $\Z/N\Z$ and 
the other on the space of bounded functions on a measure space.
We used the ergodic seminorms to define factors 
of a measure space, and then showed that these factors have 
algebraic structure.  This algebraic structure is the 
main ingredient in proving convergence of multiple ergodic averages 
along arithmetic progressions, and along other sequences.  
Gowers norms have since been used in other contexts, 
including the proof of Green and Tao~\cite{GT0} that the 
primes contain arbitrarily long arithmetic progressions.  
The connection between nilsystems in ergodic theory and the 
algebraic nature of analogous combinatorial objects 
has yet to be fully understood.  The beginning of this is carried 
out by Green and Tao
(see~\cite{GT},~\cite{GT2} and~\cite{GT3}), 
including an inverse theorem for the third Gowers norm.

In this article, we define related seminorms 
on bounded sequences and prove a structure theorem 
and an inverse theorem for it.  We also give some 
ergodic theoretic applications of these constructions.  
 These applications include 
a version of the Wiener-Wintner ergodic theorem extended 
to nilsequences, a spectral type theorem for nilsequences, 
and a weighted ergodic theorem.  
Polynomial versions of these results are contained in a 
forthcoming article.
All these properties depend on the connection to 
algebraic structures and we describe these structures more precisely.

\subsection{Nilsystems and nilsequences}

In the inverse and structure theorems described above, a key role 
is played by algebraic objects, the nilsystems:
\begin{definition}
Assume that $G$ is a $k$-step nilpotent Lie group and 
$\Gamma\subset G$ is a discrete, cocompact subgroup of $G$. 
The compact manifold $X = G/\Gamma$ is called a {\em $k$-step 
nilmanifold}. 
The \emph{Haar measure} $\mu$ of $X$ is the unique probability measure 
invariant under the action $x\mapsto g.x$ of $G$ on $X$ by left 
translations. 
Letting $T$ denote left multiplication  by the fixed 
element $\tau\in G$,  we call 
$(X, \mu, T)$ a {\em $k$-step nilmanifold}\footnote{
$X$ is endowed with its Borel $\sigma$-algebra $\CX$. In general, 
we omit the associated $\sigma$-algebra from our notation, 
writing $(X, \mu, T)$ for a measure preserving probability 
system rather than $(X, \CX, \mu, T)$.  We implicitly assume 
that all measure preserving systems are probability systems.}.  
\end{definition}

Loosely speaking, the Structure Theorem of~\cite{HK} 
states that if one wants to understand the multiple ergodic 
averages 
$$
\frac{1}{N}\sum_{n=0}^{N-1}f_1(T^nx)\ldots f_k(T^{kn}x) \ ,
$$
where $k\geq 1$ is an integer, $(X, \mu, T)$ is a measure preserving system, 
and $f_1, \ldots, f_k\in L^\infty(\mu)$, one can replace 
each function by its conditional expectation on 
some nilsystem.  Thus one 
can reduce the problem to studying the same average in a 
nilsystem, reducing averaging in an arbitrary system to a
more tractable question.  

A related problem is study of the multicorrelation sequence
$$
c_n := \int T^nf\cdot T^{2n}f\cdot\ldots\cdot T^{kn}f\,d\mu\ , 
$$
where $k\geq 1$ is an integer, 
$(X, \mu, T)$ is a measure preserving system, and 
$f\in L^\infty(\mu)$.
In~\cite{BHK}, we defined sequences that arise from nilsystems (the 
{\em nilsequences}) and show that a multicorrelation 
sequence can be decomposed into a sequence that is 
small in terms of density and 
a $k$-step nilsequence.  We define this second term precisely:

\begin{definition}
Let $(X, \mu, T)$ be a $k$-step nilsystem, $f\colon X\to\C$ a 
continuous function, $\tau\in G$, and $x_0\in X$. 
  The sequence $(f(\tau^nx_0)\colon n\in\Z)$ is a 
{\em basic $k$-step nilsequence}. 
If, in addition, the function $f$ is smooth, then the sequence 
$(f(\tau^nx_0)\colon n\in\Z)$ is called a {\em smooth $k$-step 
nilsequence}. A {\em $k$-step nilsequence} 
is a uniform limit of basic $k$-step nilsequences.

\end{definition}

The family of $k$-step nilsequences forms a closed, 
shift invariant subalgebra of 
sequences in $\ell^\infty(\Z)$.  One step nilsequences are 
exactly the almost periodic sequences.
An example of a $2$-step nilsequence is the sequence  
$(\exp(\pi in(n-1)\alpha)\colon n\in\Z)$, 
where $\alpha$ lies in the torus $\T=\R/\Z$.  
(The collection of all $2$-step nilsequences is described 
fully and classified in~\cite{HK2}.)

\subsection{Direct theorems and inverse theorems}
We define a new seminorm on bounded sequences and use this 
seminorm, an associated dual norm, and nilsequences to derive 
direct and inverse theorems.  These seminorms on $\ell^\infty(\Z)$
arise via an averaging process, and there is more than one 
natural way to take such an average. 
The first is looking along a 
particular sequence of intervals of integers whose lengths 
tend to infinity, and taking the average over 
these intervals.  This corresponds, in some sense, to a local point 
of view, as such an averaging scheme does not take into account 
what happens outside this particular sequence of intervals.  
A second way to take an average is to allow all choices of 
intervals.  This uniform point of view gives us further information 
on the original sequence.  

Averaging in $\Z$, the first version gives rise to the classic 
notion of density, taking the proportion of a set relative to the sequence
of intervals $[1, \ldots, N]$, while the second
gives rise to the slightly different
notion of Banach density, where 
the density is computed relative to any 
sequence of intervals whose lengths tend to infinity.
Each type of averaging gives rise for 
each integer $k\geq 1$
to some sort of uniformity measurement (seminorm, norm, or a 
version thereof) on bounded sequences.

We use the seminorms associate to each of these 
 averaging methods to address analogs of 
combinatorial results.
A classical problem in combinatorics is to start with a finite 
set $A$ of integers (for example) and say something about 
properties of sets that can be built from $A$, such as 
the sumset $A+A$ or product set $A\cdot A$.  Such results are 
referred to as direct theorems.
Inverse theorems start with the sumset, product set, or other 
information derived from a finite set, and then try to deduce 
information about the set itself.  

We prove both a direct theorem and an inverse theorem.  
For the direct theorem, we show that there is a dual norm 
that acts as an upper bound 
on the correlation of a bounded sequence with a nilsequence. 
We also prove an inverse theorem for the seminorms, 
showing how a bounded sequence correlates with a nilsequence.
This is an $\ell^\infty$ version of the Gowers Inverse Conjecture 
made by  Green and Tao~\cite{GT2}.
This conjecture was resolved by them for the third Gowers 
norm in~\cite{GT3}.

Using the direct theorems, we derive a weighted 
multiple ergodic convergence theorem. 
We believe that one should be able to use these methods 
to derive other combinatorial results.

The tools used in this paper have several sources.
One is a version of the Furstenberg Correspondence Principle 
(see~\cite{F}), 
used to translate the problems into ergodic theoretic statements.  
Another is the connection of the seminorms 
we define with the algebraic structure of nilsystems, 
using properties of the ergodic seminorms developed in~\cite{HK}.
Throughout, we use some harmonic analysis on nilmanifolds.

This article can be viewed as an ergodic perspective on 
the development of a ``higher order Fourier analysis'' that has 
been proposed by Green and Tao~\cite{GT2}.  Our direct 
results develop harmonic analysis relative to the 
standard Fourier analytic methods and our local 
inverse results lend support to Green-Tao conjecture 
of an inverse theorem for the Gowers norms.

\subsection{Organization of the paper}

In the next section, we define the seminorms on $\ell^\infty(\Z)$ 
and give their basic properties.
We then state the main results first for $k=2$ and then for general $k$, with 
the 
intention of clarifying the objects under study.  Section~\ref{sec:tools}
gives the background on ergodic seminorms and nilsystems.
In Section~\ref{sec:correspond}, we give a presentation 
of the Correspondence Principle that allows us to 
prove the properties of the $\ell^\infty(\Z)$ seminorms introduced 
in Section~\ref{sec:summary}.
In Section~\ref{sec:duality}, we study the dual 
norm associated to these seminorms and use it to prove 
the direct theorems on the seminorms.  We 
prove the inverse theorems in Section~\ref{sec:extension}, 
using an extension of the Correspondence Principle and 
in Section~\ref{sec:application} we give some 
ergodic theoretic consequences of these results. 
Throughout we make use of the connection with the ergodic 
seminorms.

\section{Summary of the results}
\label{sec:summary}
We introduce seminorms on $\ell ^\infty(\Z)$ corresponding to the 
Gowers norms~\cite{Gowers} in the finite setting and to the 
seminorms in ergodic theory introduced in~\cite{HK}.
We begin with some definitions and statements 
of the main properties.  
After defining the relevant seminorms, we give the statements of 
the results, beginning with the sample case of $k=2$.

\begin{notation}
We write sequences as $\ua=(a_n\colon n\in\Z)$ and 
we write the uniform norm of this sequence as $\norm\ua_\infty$.

By an {\em interval}, we mean an interval in $\Z$. If $I$ is an 
interval, $|I|$ denotes its length.

We write $z\mapsto Cz$ for complex conjugation in $\C$. Thus 
$C^kz=z$ if $k$ is an even integer and $C^kz=\bar z$ if $k$ is an odd 
integer. 

For every $k\geq 1$, points of $\Z^k$ are written
 $h=(h_1,\dots,h_k)$.
For $\epsilon = (\epsilon_1, \ldots, \epsilon_k)\in\{0,1\}^k$ 
and $h = (h_1, \ldots, h_k)\in\Z^k$, we define 
$$
|\epsilon| = \epsilon_1+\ldots+\epsilon_k\text{ and }
\epsilon\cdot h = \epsilon_1\cdot h_1+\ldots +
\epsilon_k\cdot h_k\ .  
$$
\end{notation}

Further notation on averages of sequences of intervals is
given at the end of this Section.

\subsection{The local ``seminorms'' and the 
uniformity seminorms on $\ell^\infty(\Z)$}
\label{sec:seminorms}
We define two quantities that are measurements on bounded sequences.
The proofs rely on 
material from a variety of sources (summarized in Section~\ref{sec:tools}) 
and some machinery that we develop, and so we postpone 
them until Section~\ref{sec:correspond}. 
In fact, some of the properties stated in this section can be 
proved via direct computations.   However, we prefer proofs relying 
on the Furstenberg correspondence principle, as we use a modification of this principle to prove stronger results.

We introduce the property that allows us to define certain ``seminorms.''
\begin{definition} 
\label{def:corell}
Let $k\geq 1$ be an integer, $\ua = (a_n\colon n\in\Z)$ 
be a bounded sequence, and $\uI=(I_j\colon j\geq 1)$ be 
a sequence of intervals whose lengths tend to infinity.
We say that the sequence $\ua$ 
{\em satisfies property $\CP(k)$ on $\uI$} if
for all $h=(h_1,\dots,h_{k})\in\Z^{k}$,
the limit
$$
\lim_{j\to+\infty}\frac 1{|I_j|}\sum_{n\in I_j}
\prod_{\epsilon\in\{0,1\}^{k}} C^{|\epsilon|}a_{n+ h\cdot\epsilon}
$$
exists.
We denote this limit by $c_h(\uI,\ua)$.  
\end{definition}

Given a bounded sequence $\ua$ and a sequence
of intervals whose lengths tend to 
infinity, one can always pass to a subsequence on which $\ua$ satisfies 
$\CP(k)$. 

\begin{proposition}
\label{prop:seminorms}
Let $k\geq 1$ be an integer,
$\uI=(I_j\colon j\geq 1)$ be a sequence of intervals whose lengths 
tend to infinity, and let $\ua$ be a bounded sequence satisfying 
property $\CP(k)$ on $\uI$.  
Then then  limit
$$
\lim_{H\to+\infty}\frac 1{H^k}\sum_{h_1,\dots,h_k=0}^{H-1} c_h(\uI,\ua)\ ,
$$ 
exists and is non-negative.
\end{proposition}

Using this proposition, we define:
\begin{definition}
\label{def:semi}
For  an integer $k\geq 1$, a sequence of intervals 
$\uI=(I_j\colon j\geq 1)$, and a bounded sequence $\ua$ 
satisfying  property $\CP(k)$ on $\uI$, define
$$
 \norm\ua_{\uI,k}=
\Bigl(\lim_{H\to+\infty}\frac 1{H^k}\sum_{h_1,\dots,h_k=0}^{H-1} 
c_h(\uI,\ua)\Bigr)^{1/2^k}\ .
$$
\end{definition}
We call $\norm\cdot_{\uI,k}$ a \emph{local ``seminorm''} (with quotes 
on the word seminorm), because the 
space of sequences satisfying  property $\CP(k)$ on $\uI$ is not a 
vector space.   On the other hand, we do have:

\begin{proposition}
\label{prop:seminorms2}
Assume that $k\geq 1$ is an integer, $\ua$ and $\ub$ are bounded 
sequences, and $\uI$ is a sequence of intervals whose lengths 
tend to infinity. 
If $\ua,\ub$ and $\ua+\ub$ satisfy property $\CP(k)$ on 
$\uI$, then 
$\norm{\ua+\ub}_{\uI,k}\leq\norm\ua_{\uI,k}+\norm\ub_{\uI,k}$.
\end{proposition}

The ``seminorms'' are also non-increasing with $k$:

\begin{proposition}
\label{prop:seminorms3}
If the bounded sequence $\ua$ satisfies properties $\CP(k)$ and 
$\CP(k+1)$ on the sequence of intervals 
$\uI$, then $\norm\ua_{\uI,k}\leq\norm\ua_{\uI,k+1}$.
\end{proposition}

We use the ``seminorm'' to define a measure of uniformity (a 
uniformity seminorm) on bounded sequences:

\begin{definition}
Let $\ua = (a_n\colon n\in\Z)$ 
be a bounded sequence and let $k\geq 1$ be an integer. We define 
the \emph{$k$-uniformity seminorm} $\norm\ua_{U(k)}$ to be the supremum 
of $\norm\ua_{\uI,k}$, where the supremum is taken over 
all sequences of intervals $\uI$ on which 
$\ua$ satisfies property $\CP(k)$.
\end{definition}

Using Proposition~\ref{prop:seminorms2}, by passing, if necessary, 
to subsequences of the sequences of intervals, 
we immediately deduce:
\begin{proposition}
\label{prop:uk-seminorm}
For every integer $k\geq 2$, $\norm\cdot_{U(k)}$ is a seminorm on 
$\ell^\infty(\Z)$.
\end{proposition}

\subsection{Comments on the definitions}
\subsubsection{}
The definitions of $\norm\ua_{\uI,k}$ and $\norm\ua_{U(k)}$ are very similar 
to those of the Gowers norms introduced in~\cite{Gowers} in the 
finite setting (meaning, for sequences indexed by $\Z/N\Z$).  
In the sequel, we establish analogs of 
properties of Gowers norms for the $\ell^\infty(\Z)$ seminorms.
The $\ell^\infty(\Z)$ seminorms are also close relatives of the 
ergodic seminorms of~\cite{HK}.  In the sequel we show 
that this resemblance is not merely formal; the link between 
the $\ell^\infty(\Z)$ seminorms and the ergodic seminorms 
is a basic tool of this paper.

\subsubsection{} It can be shown that in Proposition~\ref{prop:seminorms}
the averages on $[0,H-1]^k$ can be replaced by averages on any 
sequence of ``rectangles'' 
$(I_{H,1}\times\dots I_{H,k}\colon H\geq 1)$, 
where $I_{H,j}$ is an interval for every $j\in\{1, \ldots, k\}$ 
and every $H$ 
and $\min_j|I_{H,j}|\to+\infty$ as $H\to+\infty$; more 
generally we could also average over any 
F\o lner sequence in $\Z^k$.

\subsubsection{}For clarity,  we explain what the definitions mean when $k=1$. 
(We discuss $k=2$ in the next section.) 
Let $\ua=(a_n\colon n\in\Z)$ be a bounded sequence and let 
$\uI=(I_j\colon j\geq 1)$ be a sequence of intervals whose 
lengths tend to infinity.

Property $\CP(1)$ says that for every $h\in\Z$, the averages 
$$\frac 1{|I_j|}\sum_{n\in I_j}a_n\overline{a_{n+h}}$$
converge as $j\to+\infty$ and the definition of $\norm\ua_{\uI,1}$ is
$$
 \norm\ua_{\uI,1}=\Bigl(\lim_{H\to+\infty}\frac 1H\sum_{h=0}^{H-1}
\lim_{j\to+\infty}\frac 1{|I_j|}\sum_{n\in I_j}a_n\overline{a_{n+h}}
\Bigr)^{1/2}\ .
$$
Furthermore, 
$$
\norm\ua_{\uI,1}\geq\limsup_{j\to+\infty}\Bigl|\frac 
1{|I_j|}\sum_{n\in I_j}a_n\Bigr|
$$ 
and
$$
\norm\ua_{U(1)}=
\lim_{N\to+\infty}\ 
\sup_{M\in\Z}\ \Bigl| \frac 1{N} \sum_{n=M}^{M+N-1}a_n\Bigr|
\ .
$$ 
The first property follows easily from the van der Corput 
Lemma (see Appendix~\ref{ap:vdc}) and 
probably the second can also be proved directly.  
Both properties also follow from the discussion in 
Section~\ref{subsec:proofsseminorms}.

\subsubsection{}
The difference between the local ``seminorms'' and the uniformity 
seminorms is best illustrated by considering a randomly 
generated sequence.
Let $\ua=(a_n\colon n\in\Z)$ be a random sequence, where the $a_n$ 
are independent random variables, taking the values 
$+1$ and $-1$ each with probability $1/2$. 
Let $\uI=(I_j\colon j\geq 1)$ be a sequence of intervals whose lengths tend to infinity. 
Then for every integer $k$, 
the sequence $\ua$ satisfies property $\CP(k)$ 
on $\uI$ almost surely 
and $\norm\ua_{\uI,k}=0$. On the other hand, we 
have that 
$\norm\ua_{U(k)}=1$ almost surely. Indeed, for every integer $j\geq 1$ there exists an 
interval $I_j$ of length $j$ on which the sequence $\ua$ is constant 
and equal to $1$; taking $\uI$ to be this sequence of intervals, we have that
$\norm\ua_{\uI,k}=1$ for every integer $k\geq 1$. 
The apparent contradiction only arises because of the choice of uncountably many 
sequences of intervals.

\subsubsection{}
There are nontrivial bounded sequences for which the uniformity 
seminorm is $0$.  
This is illustrated by the following particular case of 
Corollary~\ref{cor:normtotergo}.
\begin{proposition}
Let $k\geq 1$ be an integer and assume that $(X,T)$ is a 
uniquely ergodic system with invariant 
measure $\mu)$ that is weakly mixing.
If $f$ is a function on $X$ with 
$\int f \,d\mu = 0$, then for every $x\in X$,
the sequence $(f(T^nx)\colon n\in\Z)$ has $0$ $k$-uniformity seminorm.
\end{proposition}

\subsection{The case $k=2$} 
\label{subsec:kequals2}
To further clarify the statements, 
we explain some of our general results in the particular case 
that $k=2$.  
These results are  
prototypes for the general case, but are simpler 
to state and prove. 
Most of these results can be proved without 
resorting to any significant machinery and we include one of the simpler proofs here.  
\begin{notation} We write $\T=\R/\Z$. For  $t\in\T$, $e(t)=\exp(2\pi 
it)$.
\end{notation}
The first result explains the role of the local ``seminorm'', namely 
that it acts as an upper bound:
\begin{proposition}
\label{prop:k2vdC}
If $\ua = (a_n\colon n\in\Z)$ is a bounded sequence satisfying
$\CP(2)$ on the sequence of intervals $\uI= (I_j\colon j\geq 1)$, then
$$
 \limsup_{j\to+\infty}\;
\sup_{t\in\T}\Bigl|\frac 1{|I_j|}\sum_{n\in I_j}a_ne(nt)\Bigr|
\leq\norm\ua_{\uI,2}\ .
$$
\end{proposition}

\begin{proof}
We can assume that $\norm\ua_\infty\leq 1$.
By the van der Corput Lemma (Appendix~\ref{ap:vdc}), Cauchy-Schwartz 
Inequality, and  another application of the van der Corput Lemma, 
we have that for all integers $j, H\geq 1$, and all $t\in\T$,
\begin{multline*}
 \Bigl|\frac 1{|I_j|}\sum_{n\in I_j}a_ne(nt)\Bigr|^4
\leq 
\Bigl(\frac {cH}{|I_j|}+\Bigl|\sum_{h=-H}^H\frac{H-|h|}{H^2}
\; \frac 1{|I_j|}\sum_{n\in I_j}a_n\overline{a_{n+h}}
\Bigr|\Bigr)^2\\
\leq 
\frac {c'H}{|I_j|} +\sum_{h=-H}^H\frac{H-|h|}{H^2}
\;\Bigl| \frac 1{|I_j|}\sum_{n\in I_j}a_n\overline{a_{n+h}}
\Bigr|^2\\
\leq \frac {c''H}{|I_j|}+\sum_{\ell=-H}^H\sum_{h=-H}^H
\frac{H-|\ell|}{H^2}\frac{H-|h|}{H^2}
\;\Bigl| \frac 1{|I_j|}\sum_{n\in I_j}a_n\overline{a_{n+h}}\,
\overline{a_{n+\ell}}a_{n+h+\ell}\Bigr|\ ,
\end{multline*}
where $c,c',c''$ are universal constants. Taking the limit as $j\to+\infty$ 
first
(recall that the sequence $\ua$ satisfies
$\CP(2)$ on the sequence of intervals $\uI$), 
and then as $H\to+\infty$, we  have the announced result.
 \end{proof}
 
We use this to show how such a sequence 
$\ua$ correlates with almost periodic sequences.  First a 
definition:

\begin{definition}
A sequence of the form $(e(nt)\colon n\in\Z)$ is called a 
\emph{complex exponential sequence}.
A sequence is a \emph{trigonometric polynomial} if it is 
a finite linear combination of complex exponential sequences. 
An \emph{almost periodic sequence} is a uniform limit of 
trigonometric polynomials.
\end{definition}

By approximation, it follows immediately from Proposition~\ref{prop:k2vdC} 
that:
\begin{corollary}
Let $\ub=(b_n\colon n\in\Z)$ be an almost periodic sequence. Then for every 
$\delta>0$, 
there exists a constant $c=c(\ub,\delta)$ such that
if a bounded sequence $\ua=(a_n\colon n\in\Z)$ 
satisfies property $\CP(2)$ on a 
sequence of intervals $\uI = (I_j\colon j\geq 1)$, then
$$
 \limsup_{j\to+\infty}\Bigl|\frac 1{|I_j|}\sum_{n\in I_j}a_nb_n\Bigr|\leq 
c\norm\ua_{\uI,2}+\delta\norm\ua_\infty\ .
$$
\end{corollary}

For some almost periodic sequences we have more precise bounds. A 
smooth almost periodic sequence $\ub = (b_n\colon n\in\Z)$ 
(that is, a smooth $1$-step nilsequence) can be written as
$$
 b_n=\sum_{m=1}^\infty \lambda_me(nt_m)\ ,
$$
where $t_m$, $m\geq 1$, are distinct elements of $\T$ 
and $\lambda_m\in \C$, $m\geq 1$, satisfy 
$$
 \sum_{m=1}^\infty |\lambda_m|<+\infty\ .
$$
We define
$$
 \nnorm\ub_2^*=\Bigl(\sum_{m=1}^\infty|\lambda_m|^{4/3}\Bigr)^{3/4}
$$
and we have that:
\begin{proposition}
\label{prop:upperkequals2}
Let $\ua=(a_n\colon n\in\Z)$ be a bounded 
sequence satisfying property $\CP(2)$ on the 
sequence of intervals $\uI = (I_j\colon j\geq 1)$ 
and $\ub=(b_n\colon n\in\Z)$ be a smooth almost 
periodic sequence.  Then, 
$$
 \limsup_{j\to+\infty}\Bigl\vert \frac{1}{|I_j|}\sum_{n\in I_j}a_nb_n
 \Bigr\vert
 \leq\norm\ua_{\uI,2}\,\nnorm\ub_2^*\ .
$$
\end{proposition}
The constant $\nnorm\ub_2^*$ here is the best possible. 
Undoubtedly, one could prove this result 
without resorting to special machinery, but we do not attempt this method as this is a particular case of a general result 
(Theorem~\ref{th:upperbound}). 
In fact we show that the norm $\nnorm\cdot_2^*$ acts as
the dual of the seminorm $\norm\cdot_{U(2)}$.

\subsection{Main results}
Let $k\geq 2$ be an integer. In section~\ref{subsec:dualnorm}, for every 
$(k-1)$-step 
nilmanifold $X$ we define a norm $\nnorm\cdot_k^*$ on the space 
$\CC^\infty(X)$ of smooth functions on $X$. 
We defer the precise definition, as it requires development of 
some further background.
Let $\ub$ be a smooth $(k-1)$-step nilsequence. Then there exists an 
ergodic $(k-1)$-step nilsystem (Corollary~\ref{cor:nilseqergodic}),
 a smooth function $f$ on $x$, and a point $x_0\in X$ 
with 
$$
b_n=f(T^nx_0)\text{ for every }n\in\Z\ .
$$
 The same sequence $\ub$ can be 
represented in this way in several manners,  with different systems, 
different starting points, and different functions, but we show 
(Corollary~\ref{cor:normsequence}) that all 
associated functions $f$ have the same 
norm $\nnorm\cdot_k^*$. Therefore we can define $\nnorm\ub_k^*=\nnorm 
f_k^*$ where $f$ is any of the possible functions.

\subsubsection{Direct results}

Using this norm, we have generalizations of the results 
already given for $k=2$:
\begin{theorem}[Direct Theorem]
\label{th:upperbound}
Let $\ua=(a_n\colon n\in\Z)$ be a bounded sequence that satisfies property $\CP(k)$ on the 
sequence of intervals $\uI = (I_j\colon j\geq 1)$. 
For all $(k-1)$-step smooth nilsequences $\ub$, we have
$$
\limsup_{j\to+\infty}\Bigl|
\frac 1{|I_j|}\sum_{n\in I_j}a_nb_n\Bigr|
\leq\norm\ua_{\uI,k}\,\nnorm \ub_k^*\ .
$$
\end{theorem}
By density, Theorem~\ref{th:upperbound} immediately implies:
\begin{corollary}
\label{cor:upperbound}
Let $\ub=(b_n\colon n\in\Z)$ be a $(k-1)$-step nilsequence
and $\delta>0$.  There exists a 
constant  $c=c(\ub,\delta)$ such that for every bounded sequence 
$\ua = (a_n\colon n\in\Z)$ 
satisfying property $\CP(k)$ on a sequence of intervals 
$\uI=(I_j\colon j\geq 1)$, we have
$$
 \limsup_{j\to+\infty}\Bigl|
\frac 1{|I_j|}\sum_{n\in I_j}a_nb_n\Bigr|
\leq c\norm\ua_{\uI,k}+\delta\norm\ua_\infty\ .
$$
\end{corollary}

Using these results, we immediately deduce uniform versions:
\begin{corollary}
\label{cor:upperbounduniform}
Let $\ub=(b_n\colon n\in\Z)$ be a smooth $(k-1)$-step nilsequence and
$\ua=(a_n\colon n\in\Z)$ be a bounded sequence.  Then
$$
 \lim_{N\to+\infty}\sup_{M\in\Z}\Bigl\vert\frac{1}{N}\sum_{n=M}^{N+M-1}
 a_nb_n\Bigr\vert
 \leq \norm\ua_{U(k)}\,\nnorm \ub_k^*\ .
$$
Let $\ub=(b_n\colon n\in \Z)$  be a $(k-1)$-step nilsequence 
and let $\delta>0$.  There exists a 
constant  $c=c(\ub,\delta)$ such that for every 
bounded sequence $\ua = (a_n\colon n\in\Z)$, 
$$
 \lim_{N\to+\infty}\sup_{M\in\Z}\Bigl\vert\frac{1}{N}\sum_{n=M}^{N+M-1}
a_nb_n\Bigr\vert
\leq c\norm\ua_{U(k)}+\delta\norm\ua_\infty\ .
$$
\end{corollary}

We refer to these results as direct results, meaning that we start with a 
sequence and derive its correlation with nilsequences.
One can view them as upper bounds, because 
they give an upper bound 
between the correlation of a sequence with a nilsequence.

\subsubsection{Inverse results}
The next results are in the opposite direction of the direct results 
of the previous section, and 
we refer to them as ``inverse results''.

\begin{theorem}[Inverse Theorem]
\label{th:inverse}
 Let $\ua = (a_n\colon n\in\Z)$ be a bounded sequence.  
 Then for every $\delta>0$, there exists a $(k-1)$-step 
smooth nilsequence $\ub=(b_n\colon n\in\Z)$  such that 
$$
 \nnorm \ub_k^*=1\text{ and }
\lim_{N\to+\infty}\,\sup_{M\in\Z}\Bigl|
\frac{1}{N}\sum_{n=M}^{M+N-1}
a_nb_n\Bigr|
\geq \norm \ua_{U(k)}-\delta\ .
$$
\end{theorem}

Summarizing this  theorem and Corollary~\ref{cor:upperbounduniform}
we have
\begin{corollary}
For every bounded sequence $\ua=(a_n\colon n\in\Z)$, 
$$
\norm\ua_{U(k)}=\sup_{\substack{\ub =(b_n)
\text{ \upshape is a smooth}\\\text{\upshape nilsequence and }
\nnorm\ub_k^*=1}} 
\lim_{N\to+\infty}\sup_{M\in\Z}\Bigl\vert\frac{1}{N}\sum_{n=M}^{N+M-1}
a_nb_n\Bigr\vert \ .
$$
\end{corollary}

This means that we can view the
norm $\nnorm\cdot_k^*$ as the dual norm of the 
uniformity seminorm $\norm\cdot_{U(k)}$.

\begin{corollary}
\label{cor:unif-zero}
For a bounded sequence $\ua=(a_n\colon n\in\Z)$, 
the following properties are equivalent:
\begin{enumerate}
\item $\norm\ua_{U(k)}=0$.
\item $\displaystyle \lim_{N\to+\infty}\sup_{M\in\Z}\Bigl\vert\frac{1}{N}\sum_{n=M}^{N+M-1}
a_nb_n\Bigl|=0$ for every $(k-1)$-step smooth nilsequence 
$\ub=(b_n\colon n\in\Z)$.
\item $\displaystyle \lim_{N\to+\infty}\sup_{M\in\Z}\Bigl\vert\frac{1}{N}\sum_{n=M}^{N+M-1}
a_nb_n\Bigl|=0$ for every $(k-1)$-step nilsequence $\ub=(b_n\colon n\in\Z)$.
\end{enumerate}
\end{corollary}

For $k=2$, Corollary~\ref{cor:unif-zero}, 
Proposition~\ref{prop:k2vdC}, and a density argument imply that the three
equivalent conditions of Corollary~\ref{cor:unif-zero} are also equivalent to
\begin{enumerate}
\setcounter{enumi}{3}
\item
For every $t\in\T$, 
$\displaystyle 
 \lim_{N\to+\infty}\,\sup_{M\in\Z}\Bigl\vert\frac{1}{N}\sum_{n=M}^{N+M-1}
a_ne(nt)\Bigl|= 0$.
\item
$\displaystyle
\lim_{N\to+\infty}\, \sup_{t\in\T}\,\sup_{M\in\Z}\Bigl|
\frac 1N\sum_{n=M}^{M+N-1}a_ne(nt)\Bigr|= 0$.
\end{enumerate}

\subsubsection{A counterexample}
It is important to note that the inverse results have no version 
involving local ``seminorms'' and we give here an example  
illustrating this point.

Let $(N_j\colon j\geq 1)$ be an increasing sequence of integers with $N_1=0$ and 
tending sufficiently fast to 
$+\infty$.
For $j\geq 1$ let $I_j=[N_j,N_{j+1}-1]$  and  let $\uI=(I_j\colon j\geq 
1)$.
 Let the sequence $\ua$ be defined by $a_n=e(n/j)$ if 
$N_j\leq|n|<N_{j+1}$.
Then  $\norm\ua_{\uI,2}=1$ and for every 
$t\in\T$, the average of $a_ne(nt)$ on the interval $I_j$ converges to 
zero as $j\to+\infty$.
Therefore, for every almost periodic sequence $\ub$, the average of 
$a_nb_n$ on $I_j$ also converges to zero.

This highlights a difference between the finite case, where the 
norms are defined on $\Z/N\Z$, and the infinite case. 
One can not construct such a sequence where the behavior worsens
as one tends to infinity. 

\subsection{A condition for convergence}
\label{subsec:condconverge}
\begin{theorem}
\label{th:cond-conv}
For a bounded sequence $\ua= (a_n\colon n\in\Z)$, 
the following are equivalent.
\begin{enumerate}
\item
For every $\delta >0$, the sequence $\ua$ can be written as
$\ua'+\ua''$ where $\ua'$ is a $(k-1)$-step nilsequence
and $\norm{\ua''}_{U(k)}<\delta$.
\item
For every $(k-1)$-step nilsequence $\uc=(c_n\colon n\in\Z)$, 
the averages of $a_nc_n$
converge, meaning that the limit
$$
\lim_{j\to+\infty}\frac 1{|I_j|}\sum_{n\in I_j} a_nc_n
$$
exists for every sequence $(I_j\colon j\geq 1)$ of intervals whose 
lengths tend to infinity.
\end{enumerate}
\end{theorem}

In Proposition~\ref{prop:riemann2}, we give
a method to build
sequences satisfying the (equivalent) properties of
Theorem~\ref{th:cond-conv},
checking that the sequences verify the first property.  
As this proposition uses material not yet defined, we do not 
state it here but only give two examples of its application.

A generalized polynomial is defined to be a real valued function 
that  is obtained from the identity function and real constants by 
using (in arbitrary order) the operations of 
addition, multiplication, and taking the integer part. We have:
\begin{proposition}
\label{prop:genpol}
Let $p$ be a generalized polynomial and for every $n\in\Z$, let 
$\{p(n)\}$ be the fractional part of $p(n)$. Then the sequences 
$(\{p(n)\}\colon n\in\Z)$ and $(e(p(n))\colon n\in\Z)$ 
satisfy the (equivalent) properties of Theorem~\ref{th:cond-conv}.
\end{proposition}

The \emph{Thue-Morse sequence} $\ua = (a_n\colon n\in\Z)$ 
is given by $a_n=1$ if the sum of the 
digits of $|n|$ written in base $2$ is odd and $a_n=0$ otherwise.
In Section~\ref{subsec:morse} we show:
\begin{proposition}
\label{prop:morse}
The Thue-Morse sequence satisfies the properties of 
Theorem~\ref{th:cond-conv}.
\end{proposition}
A similar method can be used for other 
sequences, for example for all sequences associated to primitive substitutions 
of constant length (see~\cite{Q} for the definition).

\subsection{An application to ergodic theory}
\subsubsection{}We recall a classical result in ergodic theory.

\begin{theorem*} [Wiener-Wintner ergodic theorem~\cite{WW}] 
Let  $(X,\mu,T)$ be an ergodic system and $\phi\in 
L^\infty(\mu)$. Then
there exists $X_0\subset X$ with $\mu(X_0)=1$ 
such that 
$$
 \frac 1N\sum_{n=0}^{N-1} \phi(T^nx)\,e(nt)
$$
converges for every $x\in X_0$ and every $t\in\T$.
\end{theorem*}
The important point here is that the set $X_0$ does not depend on 
the choice of $t$. 

We also recall an immediate corollary of the spectral theorem:
\begin{corollary*}[of the spectral theorem]
Let  $\ua=(a_n\colon n\in\Z)$ be a bounded 
sequence and assume that 
$$
 \lim_{N\to+\infty}\sum_{n=0}^{N-1} a_ne(nt)
$$
exists for every $t\in\T$. Then for every system $(Y,\nu,S)$ and 
every $f\in L^2(\nu)$, the averages
$$
 \frac 1N\sum_{n=0}^{N-1}a_nS^n f 
$$
converge in $L^2(\nu)$ as $N\to+\infty$.
\end{corollary*}

Putting these two results together, we have:
\begin{corollary*} Assume that $(X,\mu,T)$ is an ergodic system and 
$\phi\in L^\infty(\mu)$. There exists $X_0\subset X$ with 
$\mu(X_0)=1$ such that for every $x\in X_0$, 
every system $(Y,\nu,S)$, and every $f\in L^2(\nu)$, the averages
$$
 \frac 1N\sum_{n=0}^{N-1}\phi(T^nx)S^n f 
$$
converge in $L^2(\mu)$ as $N\to+\infty$.
\end{corollary*}

The strength of this result is 
that the set $X_0$ does not depend either on $Y$ or on $f$.  
We say that for every 
$x\in X_ 0$, the sequence $(\phi(T^nx))$ is 
a \emph{universally good for the convergence in mean of 
ergodic averages}.
In fact, for almost every $x$, this sequence is also universally good 
for the almost everywhere convergence~\cite{BFKO}, but 
we do not address this strengthening here.  

\subsubsection{} We generalize these results for multiple ergodic averages. 
 We start with a 
generalization of the Wiener-Wintner Theorem, where 
we can replace the exponential sequence $e(nt)$ by 
an arbitrary nilsequence. 

\begin{theorem}[A generalized Wiener-Wintner Theorem]
\label{th:NilWW}
Let  $(X,\mu,T)$ be an ergodic system and $\phi$ be a bounded measurable 
function on $X$. 
Then there exists $X_0\subset X$ with $\mu(X_0)=1$ such that
for every $x\in X_0$, the averages 
$$
 \frac 1N\sum_{n=0}^{N-1}\phi(T^nx)\,b_n 
$$
converge as $N\to+\infty$ 
for every $x\in X_0$ and every nilsequence $\ub=(b_n\colon n\in\Z)$.  
\end{theorem}

We give a sample application. Generalized polynomials were 
defined in Section~\ref{subsec:condconverge}.

\begin{corollary}
\label{cor:WW} 
Let  $(X,\mu,T)$ be an ergodic system, $\phi$ be a bounded measurable 
function on $X$, and $X_0$ be the subset of $X$ introduced in 
Theorem~\ref{th:NilWW}. Then for every $x\in X_0$ and every 
generalized polynomial $p$,  the averages 
$$
 \frac 1N\sum_{n=0}^{N-1}\phi(T^nx)\{p(n)\}\text{ and }
 \frac 1N\sum_{n=0}^{N-1}\phi(T^nx)e(p(n))
$$
converge.
\end{corollary}
(Recall that $\{p(n)\}$ denotes the fractional part of $p(n)$.)
For standard polynomial sequences, this result was proven 
by Lesigne~\cite{lesigne2}.

We next have a version of the spectral result for 
higher order nilsequences:

\begin{theorem}[A substitute for the corollary of the Spectral Theorem]
\label{th:convnilconvmult}
Let $k\geq 1$ be an integer and $\ua=(a_n\colon n\in\Z)$ be a bounded 
sequence such that the averages
$$
 \frac 1N\sum_{n=0}^{N-1}a_nb_n
$$
converge as $N\to+\infty$ for every $k$-step 
nilsequence $\ub=(b_n\colon n\in\Z)$. 
Then for every system $(Y,\nu,S)$ and every
 $f_1,\dots,f_k\in L^\infty(\nu)$, the averages
\begin{equation}
\label{eq:nilconvmult}
\frac 1N\sum_{n=0}^{N-1}a_n\,S^nf_1.S^{2n}f_2.\cdots.S^{kn}f_k
\end{equation}
converge in $L^2(\nu)$.  
\end{theorem}
Combining these theorems, we immediately deduce:
\begin{theorem}
\label{th:multipleuniversal}
Let $(X,\mu,T)$ be an ergodic system and $\phi\in L^\infty(\mu)$.
Then there exists $X_0\subset X$ with $\mu(X_0)=1$ such that
for every $x_0\in X$, every system $(Y,\nu,S)$, every integer $k\geq 
1$, and all functions $ f_1,\dots,f_k\in L^\infty(\nu)$, the 
averages
$$
\frac 1N\sum_{n=0}^{N-1}\phi(T^nx)\,S^nf_1.S^{2n}f_2.\cdots.S^{kn}f_k
$$ 
converge in $L^2(\nu)$ as $N\to+\infty$.
\end{theorem}
In short, for  every 
$x\in X_0$, the sequence $(\phi(T^nx))$ is 
\emph{universally good for the convergence in mean of multiple 
ergodic averages}.

While Theorems~\ref{th:NilWW} and~\ref{th:convnilconvmult} are 
results about nilsequences, nilsequences do not 
appear in the statement of Theorem~\ref{th:multipleuniversal}:  they 
occur only as tools in the proof, playing the role of complex exponentials 
in the classical results.

By successively using Theorems~\ref{th:cond-conv} 
and~\ref{th:convnilconvmult}, we obtain further examples of 
universally good sequences for the convergence in mean of multiple 
ergodic averages. For example, by 
Proposition~\ref{prop:genpol}, for every generalized polynomial 
$p$ the sequence $(\{p(n)\}\colon n\in\Z)$ is a 
universally good sequence for the convergence in mean of multiple 
ergodic averages, as 
is the sequence $(e(p(n))\colon n\in\Z)$. 
By Proposition~\ref{prop:morse}, so is the
Thue-Morse sequence.

\subsection{Some notation for averages}
In this paper we continuously take limits of averages on sequences of 
intervals. Writing the cumbersome formulas or replacing 
them by long explanations would make the paper unreadable and so 
we introduce 
some short notation. However, we 
continue using explicit formulas in the main statements.

We have several different notions of averaging
for a sequence in $\ell^\infty(\Z)$: over a particular
sequence of intervals or uniformly over all intervals.

For averaging over a particular sequence of intervals, we define:
\begin{definition}
Let $\ua=(a_n\colon n\in\Z)$ be a bounded sequence and
let $\uI=(I_j\colon j\geq 1)$ be a sequence of intervals 
whose lengths $|I_j|$ tend to infinity.
Define
$$
 \LimsupAvI{\uI}{a_n}=\limsup_{j\to+\infty}\Bigl|\frac
1{|I_j|}\sum_{n\in I_j}a_n\Bigr|\ .
$$
The \emph{averages of the sequence $\ua$ on $\uI$ converge} if the limit
$$
 \lim_{j\to+\infty}\frac 1{|I_j|}\sum_{n\in I_j}a_n
$$
exists.  We denote this limit by $\LimAvI{\uI}{a_n}$ and call
this the \emph{average over $\uI$ of the sequence $\ua$}.
\end{definition}

For taking a uniform average, we define:
\begin{definition}
Let $\ua=(a_n\colon n\in\Z)$ be a bounded sequence.
The \emph{upper limit of the averages of the sequence $\ua$} is defined to be
$$
\LimsupAv{a_n}=\lim_{N\to+\infty}\
\sup_{M\in\Z}\ \Bigl| \frac 1{N} \sum_{n=M}^{M+N-1}a_n\Bigr|\ .
$$
(Note that this limit exists by subadditivity.)

The \emph{averages of the sequence $\ua$ converge} if
the limit $\LimAvI{\uI}{a_n}$ exists for all
sequences of intervals $\uI=(I_j\colon j\geq 1)$
whose lengths $|I_j|$ tend to infinity.
We denote this (common) limit by $\LimAv {a_n}$
and call this the \emph{uniform average of the sequence $\ua$}. 
\end{definition}

Assuming the existence of the uniform average, it follows that 
$$
 \lim_{N\to+\infty}\
\sup_{M\in\Z}\ \Bigl|\LimAv{a_n}- \frac 1{N} 
\sum_{n=M}^{M+N-1}a_n\Bigr|=0\ .
$$

\section{Some tools}
\label{sec:tools}

\subsection{Nilmanifolds and nilsystems}
\subsubsection{The definitions}
Short definitions were given in the introduction and we repeat them 
here in a more complete form.

Let $G$ be a group. 
For $g,h\in G$, we write $[g,h]=ghg\inv h\inv$ for the 
commutator of $g$ and $h$ and we write $[A,B]$ for the 
subgroup spanned by $\{[a,b]\colon a\in A, b\in B\}$.
The {\em commutator subgroups} $G_j$, $j\geq 1$, 
are defined inductively by setting 
$G_1=G$ and $G_{j+1}=[G_j,G]$.  
Let $k\geq 1$ be an integer. We say that $G$ is \emph{$k$-step 
nilpotent} if $G_{k+1}$ is the trivial subgroup.

Let $G$ be a $k$-step nilpotent Lie group and $\Gamma$ a discrete 
cocompact subgroup of $G$.  
The compact manifold $X=G/\Gamma$ is 
called a \emph{$k$-step nilmanifold}.  
The group $G$ acts on $X$ by 
left translations and we write this action as 
$(g,x)\mapsto g.x$.  The 
{\em Haar measure} $\mu$ of $X$ is the unique probability measure on 
$X$ invariant under this action.

Let $\tau\in G$ and $T$ be the transformation $x\mapsto \tau.x$ of 
$X$. Then $(X,T,\mu)$ is called a {\em $k$-step nilsystem}.  
When the measure is not needed for results, we omit and write that 
$(X,T)$ is a $k$-step nilsystem.

Nilsystems are \emph{distal} topological dynamical systems. This 
means that, if $d_X$ is a distance on $X$ defining its topology,
then for every $x,x'\in X$,
$$
 \text{if } x\neq x', \text{ then }\inf_{n\in\Z}d_X(T^ny,T^ny')>0\ .
$$

Let $f$ be a continuous (respectively, smooth) function on $X$ and $x_0\in 
X$. The sequence $(f(T^nx_0)\colon n\in\Z)$ is called a 
\emph{basic} (respectively, \emph{smooth}) \emph{$k$-step nilsequence}.
A $k$-step nilsequence is a uniform limit of basic $k$-step 
nilsequences. Therefore, smooth $k$-step nilsequences are dense in the 
space of all $k$-step nilsequences under the uniform norm.

The Cartesian product of two $k$-step nilsystems is again a $k$-step 
nilsystem. It follows that the space of $k$-step nilsequences is an 
algebra under pointwise addition and multiplication. Moreover, this 
algebra is invariant under the shift. 

As an example, $1$-step nilsystems are translations on compact 
abelian Lie groups and $1$-step nilsequences are exactly almost periodic 
sequences. For examples of $2$-step nilsystems and a detailed study 
of $2$-step nilsequences, 
see~\cite{HK2}. 

A general reference on nilsystems is~\cite{AGH} and the results 
summarized in the next few sections are contained in the literature.  See, 
for example~\cite{lesigne} and~\cite{leibman}.

\subsubsection{Ergodicity}
\label{subsec:ergonil}
\begin{theorem}
\label{th:ergonil} Let $k\geq 1$ be an integer. 
For a $k$-step nilsystem $(X=G/\Gamma, T)$ with Haar measure $\mu$, 
the following properties are equivalent: 
\begin{enumerate}
\item
$(X,T)$ is transitive, meaning that it admits a dense orbit.
\item
$(X,T)$ is minimal, meaning that every orbit is dense.
\item
$(X,T)$ is uniquely ergodic.
\item
$(X,\mu, T)$ is ergodic.
\end{enumerate}
\end{theorem}
When these properties are satisfied, we say that the system is 
\emph{ergodic}, 
even in statements of topological nature (that is, without mention 
of the measure).

\begin{theorem}
\label{th:orbitergo}
Let $k\geq 1$ be an integer,  $(X=G/\Gamma, T)$ be a $k$-step 
nilsystem where $T$ is the translation by 
$\tau\in G$. Let  $x_0\in X$ and let $Y$ be the closed orbit of $x_0$, 
meaning that $Y$ is the closure of the orbit $\{T^n x_0\colon n\in\Z\}$.
Then $(Y,T)$ is a $k$-step nilsystem. More precisely, there exist a 
closed subgroup $G'$ of $G$ containing $\tau$, such that 
$\Gamma'=\Gamma\cap G'$ is cocompact in $G'$ and $Y=G'/\Gamma'$.
\end{theorem}

If $(f(T^nx_0)\colon n\in\Z)$ is a basic (respectively, smooth) 
nilsequence, by substituting the closed orbit of $x_0$ for $X$, we 
deduce:
\begin{corollary}
\label{cor:nilseqergodic}
For every basic (respectively, smooth) $k$-step nilsequence $\ua = 
(a_n\colon n\in\Z)$, there 
exists an ergodic $k$-step nilsystem $(X,T)$, $x_0\in X$, and a 
continuous (respectively, smooth) function $f$ on $X$ with $a_n=f(T^nx_0)$ 
for every $n\in\Z$.
\end{corollary}

\begin{corollary}
\label{cor:nilseaverages}
Let $\ua=(a_n\colon n\in\Z)$ be a nilsequence. Then the averages of 
$\ua$ converge. 
\end{corollary}
\begin{proof}
By density, we can restrict to the case that $\ua$ is a basic
nilsequence, and we write it as in Corollary~\ref{cor:nilseqergodic}.
By unique ergodicity of $(X,T)$, the averages converge to $\int f\,d\mu$, 
where $\mu$ is the Haar measure of $X$.
\end{proof}

\subsubsection{A criteria for ergodicity}
\begin{theorem}
\label{th:nilergo}
 Let $k\geq 1$ be an integer, $(X=G/\Gamma,T)$ be a $k$-step
nilsystem, and assume that $T$ is 
translation by  $\tau\in G$. Assume that
\begin{itemize}
\item[(*)] \rm The group $G$ is spanned by the connected component $G_0$ 
of its unit and by $\tau$.
\end{itemize}
Then $(X,T)$ is ergodic if and only if the translation induced by 
$\tau$ on the compact abelian group $Z=G/G_2\Gamma$ is ergodic.
\end{theorem}

Conversely, let $(X=G/\Gamma,T)$ be an ergodic nilsystem where $T$ is 
the translation by  $\tau\in G$. Let $G_1$ be the subgroup spanned by 
$G_0$ and $\tau$ and set $\Gamma_1=\Gamma\cap G_1$. Then $G_1$ is an open 
subgroup of $G$, $\Gamma_1$ is a discrete cocompact subgroup of $G_1$, 
and by ergodicity,  the image of $G_1$ in $X$ under 
the natural projection is onto. We can therefore identify $X$ with 
$G_1/\Gamma_1$.  Thus we can assume that hypothesis (*) of 
Theorem~\ref{th:nilergo} is satisfied. 
Throughout this paper, we implicitly assume that this hypothesis holds.

\subsubsection{The case of several commuting transformations}
Let $X=G/\Gamma$ be a nilmanifold and  let $\tau_1,\dots,\tau_\ell$ be 
commuting elements of $G$. For $1\leq i\leq \ell$ let $T_i\colon 
X\to X$ be the translation by $\tau_i$. Then the results of 
Section~\ref{subsec:ergonil} still hold, modulo the obvious changes.  
We do not give the modified statements here, with the exception of 
Theorem~\ref{th:nilergo}:
\begin{theorem}
\label{th:multiergo}
Let $X=G/\Gamma$ be a nilmanifold, $\tau_1,\dots,\tau_\ell$ be 
commuting elements of $G$, and for $1\leq i\leq \ell$ let $T_i\colon 
X\to X$ be the translation by $\tau_i$. Assume that:
\begin{itemize}
\item[(**)] \rm The group $G$ is spanned by the connected component $G_0$ 
of its unit and by $\tau_1,\dots,\tau_\ell$.
\end{itemize}
Then $X$ is ergodic under the action of $T_1,T_2,\dots,T_\ell$ 
if and only if the action induced by these transformations 
on the compact abelian group $Z=G/G_2\Gamma$ is ergodic.
\end{theorem}

\subsection{The measures $\mu\type k$ and HK-seminorms}
\label{sec:ergseminorms}\strut

In the rest of this section we consider arbitrary ergodic systems
and we assume that $k\geq 1$ is an integer. 
We review the construction and properties 
of certain objects on $X^{2^k}$ defined in~\cite{HK}.

\subsubsection{Some notation}
We introduce some notation to keep track of the $2^k$ copies of $X$.
If $X$ is a set, we write $X\type k = X^{2^k}$ and index 
these copies of $X$ by $\{0,1\}^k$.
An element of $X\type k$ is written as 
$$
\ulx=(x_\epsilon\colon\epsilon\in\{0,1\}^k) \ .
$$
We recall that for $\epsilon\in\{0,1\}^k$ and 
$h\in\Z^k$, we write $|\epsilon|=\epsilon_1+\dots+\epsilon_k$ and 
$\epsilon\cdot h=\epsilon_1h_1+\dots+\epsilon_kh_k$.

We write the element with all $0$'s of $\{0,1\}^k$ as
$\bzero=(0,0,\dots,0)$. 
We often give the $\bzero$-th coordinate of a point of
 $X\type k$ a 
distinguished role and we 
write 
$$
 X\type k=X\times X\type k_*\  ,\text{ where } X\type k_*=X^{2^k-1}\ .
$$
The coordinates of $X\type k_*$ are 
 indexed by the set 
$$
\{0,1\}_*^k = \{0,1\}^k\setminus\{\bzero\}
$$
and  a point of $X\type k$ is often written
$$
 \ulx = (x_\bzero, \ulx_*)\ , \text{ where }
\ulx_* = 
(x_\epsilon\colon \epsilon\in\{0,1\}^k_*)\ .
$$

When $(X,\mu, T)$ is a measure preserving system, 
we also have notation for some transformations 
that are naturally defined on $X\type k$.  Namely, 
we write $T\type k$ for the transformation $T\times T\times\ldots\times T$, 
taken $2^k$ times. Moreover, if $i\in\{1, \ldots, k\}$, we define 
$$
({T_i}\type k\ulx)_\epsilon = \begin{cases}\displaystyle
T(x_\epsilon) &\text{ if }\ \epsilon_i=1\\
  x_\epsilon &\text{ otherwise } \ .\end{cases} 
$$
For convenience, we also write $X\type 0=X$ and 
$T\type 0=T$.

\subsubsection{Measures and  HK-seminorms}
Throughout the rest of this section, $(X,\mu,T)$ denotes an ergodic system.

By induction, for every integer $k\geq 0$ we define
a measure $\mu\type k$ on $X\type k$ that is invariant 
under $T\type k$. 
We set $\mu\type 0 = \mu$. 
For $k\geq 1$,  making the natural identification of $X\type{k}$ 
with $X\type {k-1}\times X\type {k-1}$,
we write $\ulx = (\ulx', \ulx'')$ for a point of 
$X\type {k}$, with $\ulx', \ulx''\in X\type {k-1}$.  
Let $\CI\type {k-1}$ denote 
the invariant $\sigma$-algebra of the system 
$(X\type k,\mu\type {k-1}, T\type {k-1})$. 
We define $\mu\type{k}$ to be the {\em relatively independent 
joining} of $\mu\type {k-1}$ with itself over $\CI\type {k-1}$, 
meaning that if $F,G$ are bounded functions on $X\type {k-1}$, 
then
$$
\int_{X\type{k}}F(\ulx')G(\ulx'')\,d\mu\type{k}(\ulx) 
= \int_{X\type {k-1}}\E(F\mid\CI\type 
{k-1})(\uly)\cdot\E(G\mid\CI\type {k-1})(\uly)
\,d\mu\type {k-1}(\uly) \ .
$$

By induction, all the marginals of $\mu\type k$ (that is, the images 
of this measure under the natural projections $X\type k\to X$) are 
equal to $\mu$.

Since $(X\type 0,\mu\type 0, T\type 0)=(X,\mu,T)$ is ergodic, $\CI\type 
0$ is the trivial $\sigma$-algebra and $\mu\type 1=\mu\times\mu$. But 
for $k\geq 2$ the system $(X\type{k-1} ,\mu\type {k-1}, T\type 
{k-1})$ is not necessarily ergodic and $\mu\type k$ is not in general 
the product measure.

For $k\geq 1$ and every $f\in L^\infty(\mu)$,
\begin{multline*}
 \int_{X\type{k}}\prod_{\epsilon\in\{0,1\}^{k}}C^{|\epsilon|}f(x_\epsilon)
\,d\mu\type 
{k}(\ulx) \\
= \int_{X\type{k-1}}\Bigl|\E\Bigl(\prod_{\eta\in\{0,1\}^{k-1}}C^{|\eta|}
f(y_\eta)\Big\vert\CI\type{k-1}\Bigr)\Bigr|^2\,d\mu\type {k-1}(\uly) \geq 0
\end{multline*}
and so we can define the {\em HK-seminorm}
$$
 \nnorm f_k=\Bigl( \int_{X\type{k}}\prod_{\epsilon\in\{0,1\}^{k}}C^{|\epsilon|}f(x_\epsilon)\,d\mu\type 
{k}(\ulx)\Bigr)^{1/2^k}\ .
$$
To avoid ambiguities when several measures are present,
we sometimes write $\nnorm f_{\mu,k}$ instead 
of $\nnorm f_k$.

In~\cite{HK}, we 
show that $\nnorm\cdot_k$ is a seminorm on $L^\infty(\mu)$.
These seminorms satisfy an inequality similar to  the Cauchy-Schwartz-Gowers 
inequality for Gowers norms.  Namely, let 
$f_\epsilon$, $\epsilon\in\{0,1\}^k$, be $2^k$ bounded functions 
on $X$. Then
\begin{equation}
\label{eq:HKCSG}
\Bigl|\int\prod_{\epsilon\in\{0,1\}^k}f_\epsilon(x_\epsilon)\,d\mu\type k(\ulx)\Bigr|
\leq\prod_{\epsilon\in\{0,1\}^k}\nnorm {f_\epsilon}_k\ .
\end{equation}

We also have that consecutive HK-seminorms satisfy $\nnorm 
f_{k+1}\geq\nnorm f_k$, and by an application of the ergodic theorem, 
\begin{equation}
\label{eq:k-to-kplus1}
\nnorm f_{k+1}=\lim_{H\to+\infty}\Bigl(\frac 1H \sum_{h=0}^{H-1}\nnorm{ 
T^hf\cdot f}_k^{2^k}\Bigr)^{1/2^{k+1}}\ .
\end{equation}

Using the definition and the fact that the marginals of 
$\mu\type k$ are equal to $\mu$, we have that for all $f\in L^{2^k}(\mu)$, 
\begin{equation}
\label{eq:boundnnorm}
\nnorm f_k\leq\norm f_{L^{2^k}(\mu)}\ .
\end{equation}
In fact, the definition of the seminorm $\nnorm\cdot_k$ can be 
extended to $L^{2^k}(\mu)$ with the same properties.

\subsection{Convergence results}
\subsubsection{Averaging along parallelepipeds}
These seminorms and a geometric description of the factors 
they define are used to show:
\begin{theorem}[\cite{HK}, Theorem 13.1]
\label{th:cubes}
Let $f_\epsilon$, $\epsilon \in\{0,1\}^k_*$ be $2^k-1$ functions in 
$L^\infty(\mu)$. Then the averages 
$$
\frac 1{H^k}\sum_{h_1,\dots,h_k=0}^{H-1} \prod_{\epsilon\in\{0,1\}^k_*} T^{\epsilon\cdot h}f_\epsilon
$$
converge in $L^2(\mu)$ and the limit $g$ of these averages 
is characterized by 
$$
 \int h\, g\, d\mu=\int h(x_\bzero)
\prod_{\epsilon\in\{0,1\}^k_*}
 f_\epsilon(x_\epsilon)\,d\mu\type k(\ulx)
$$
for every $h\in L^\infty(\mu)$.
\end{theorem} 

In fact, we could replace the averages on $[0,H-1]^k$ by averages over 
any F\o lner sequence in $\Z^k$.
Applying Theorem~\ref{th:cubes} to the case that 
$f_\epsilon=C^{|\epsilon|}f$ for every $\epsilon$, we obtain:

\begin{corollary}
\label{cor:dual}
For every $f\in L^\infty(\mu)$, the averages
\begin{equation}
\label{eq:defDk}
 \frac 1{H^k}\sum _{h_1,\dots,h_k=0}^{H-1}
\prod_{\epsilon\in\{0,1\}^k_*} C^{|\epsilon|} f(T^{\epsilon\cdot 
h}x)
\end{equation}
converge in $L^2(\mu)$ as $H\to+\infty$.  
\end{corollary}
This leads us to a definition:
\begin{definition}
We denote the limit of~\eqref{eq:defDk} by 
$\CD_kf$ and call this function the {\em dual function of }$f$.
\end{definition}

It follows that the dual function $\CD_kf$ satisfies:
\begin{equation}
\label{eq:charactDk}
 \int \CD_k f \,.\, h \,d\mu=\int h(x_\bzero)
 \prod_{\epsilon\in\{0,1\}^k_*} C^{|\epsilon|} 
f(x_\epsilon)\,d\mu\type k(\underline x)
\end{equation}
for every $h\in L^\infty(\mu)$.  
In particular,  we have 
\begin{equation}
\label{eq:charactDk2}
\nnorm f_k^{2^k}=  \int \CD_k f\,.\,f\, d\mu=
\lim_{H\to+\infty}\frac
1{H^k}\sum_{h_1,\dots,h_k=0}^{H-1}\int\prod_{\epsilon\in\{0,1\}^k}
C^{|\epsilon|}T^{\epsilon\cdot h}f\,d\mu
\ .
\end{equation}

The notion of a dual function is implicit in~\cite{HK} 
and this notation is not used there. 
However, the notation is coherent with 
that used in several papers of Green and Tao,
where similar functions (in 
the finite setting) are called dual functions. 

The definition extends to functions in $L^{2^k}(\mu)$, for which we
use the same notation.
Indeed,  by~\eqref{eq:HKCSG},~\eqref{eq:boundnnorm}, and
density,  for $f\in L^{2^k}(\mu)$ the
convergence~\eqref{eq:defDk} holds in $L^{2^k/(2^k-1)}(\mu)$; the
limit function $\CD_kf$ belongs to
$L^{2^k/(2^k-1)}(\mu)$ with
$$
\norm{\CD_kf}_{L^{2^k/(2^k-1)}(\mu)}\leq\norm f_{L^{2^k}(\mu)}^{2^k-1}
$$
and
formula~\eqref{eq:charactDk} holds for every $h\in L^{2^k}(\mu)$.
Moreover,  $\CD_k$ is a continuous map from 
$L^{2^k}(\mu)$ to $L^{2^k/(2^k-1)}(\mu)$.

\subsubsection{Application to sequences}
\label{subsubsec:apssequences}

Let $f$ be a bounded function on $X$.  We
consider the quantities associated
to the bounded sequence $(f(T^nx)\colon 
n\in\Z)$ for a generic point $x$ of $X$, 
as in Section~\ref{sec:seminorms}.
>From  the definition of the ergodic
seminorms, the pointwise ergodic theorem,
and~\eqref{eq:charactDk2}, we immediately deduce:
\begin{corollary}
\label{cor:normergo}
Let $k\geq 2$ be an integer and let 
$\uI$ be the sequence of intervals $([0,N-1]\colon N\geq 1)$. 
Let $(X,\mu,T)$ be an ergodic system and let $f\in L^\infty(\mu)$. 
Then for almost every $x\in X$, the sequence $(f(T^nx)\colon n\in\Z)$ 
satisfies property $\CP(k)$ on $\uI$ and 
\begin{equation}
\label{eq:normergo}\norm{(f(T^nx)\colon n\in\Z)}_{\uI,k}=\nnorm f_k\ .
\end{equation}
\end{corollary}
\begin{corollary}
\label{cor:normtotergo}
Let $k\geq 2$ be an integer, 
let $(X,T)$ be a uniquely ergodic system with invariant measure $\mu$, 
and let $f$ be a Riemann integrable function on $X$.
Then for every $x\in X$ and every sequence of intervals $\uI$ whose lengths tend to infinity, 
the  sequence $(f(T^nx)\colon n\in\Z)$ satisfies property 
$\CP(k)$ on $\uI$ and  
equality~\eqref{eq:normergo} holds.

In particular, for every $x\in X$, 
$$
 \norm{(f(T^nx)\colon n\in\Z)}_{U(k)}=\nnorm f_k\ .
$$
\end{corollary}

\begin{proof}
The hypothesis means that for every $\delta>0$ there 
exists two continuous functions $g,g'$ on $X$ with $g\leq f\leq g'$ 
and $\int(g'-g)\,d\mu<\delta$. This implies that for every $h\in\Z^k$ 
the function in the last integral of formula~\eqref{eq:charactDk2} is 
also Riemann integrable. Therefore the ergodic averages of this 
function converge everywhere to its integral.
\end{proof}

\subsection{The structure Theorem}
We use the following version of the Structure Theorem of~\cite{HK}, 
which is a combination of statements in Lemma 4.3, Definition 4.10 
and Theorem 10.1 of that paper.

\begin{theorem*}[Structure Theorem]
\label{th:structure}
Let $(X,\mu,T)$ be an ergodic system. Then for every $k\geq 2$ there 
exists a system $(Z_k, \mu_k,T)$ and a factor map 
$\pi_k\colon X\to Z_k$ with the following properties:
\begin{enumerate}
\item
$(Z_k, \mu_k,T)$ is the inverse limit of a sequence of $(k-1)$-step 
nilsystems.
\item
For every function $f\in L^\infty(\mu)$, $\nnorm{f-\E(f\mid 
Z_k)\circ\pi_k}_k=0$\ .
\end{enumerate}
\end{theorem*}

Since $\nnorm f_{k+1}\geq\nnorm f_k$ for every $f\in 
L^\infty(\mu)$, the factors $Z_k$ are nested: $Z_k$ is a factor of 
$Z_{k+1}$.

We use this theorem via the following immediate corollary.
\begin{corollary}
\label{cor:structure}
Let  $(X,\mu,T)$ be an ergodic system and $f\in L^\infty(\mu)$. Then 
for every $\delta >0$, there exists a $(k-1)$-step ergodic nilsystem 
$(Y,S,\nu)$, a (measure theoretic) factor map $p\colon X\to Y$, and 
a continuous function $h$ on $Y$ with $\nnorm{f-h\circ p}_{\mu,k}<\delta$.
\end{corollary}

\section{The correspondence principle and the ``seminorms''}
\label{sec:correspond}

\subsection{The classic Correspondence Principle}
\label{sec:correspondence}

In translating Szemer\'edi's Theorem into a problem 
in ergodic theory, Furstenberg introduced the Correspondence 
Principle in~\cite{F}.  
We give a not completely classical presentation of this 
principle, which is amenable to modification in the sequel.

By a \emph{separable subalgebra of $\ell^\infty(\Z)$}, we mean a 
unitary subalgebra of $\ell^\infty(\Z)$, invariant under the shift 
 and under complex conjugation, closed in $\ell^\infty(\Z)$ and 
separable for the uniform norm written $\nnorm\cdot_\infty$.
In the sequel, we mostly consider the case of the separable subalgebra 
$\CA(\ua)$ spanned by a bounded sequence $\ua=(a_n\colon n\in\Z)$.

We write $\sigma$ for the shift on $\ell^\infty(\Z)$, and 
thus for a sequence $\ua=(a_n\colon n\in\Z)$, $\sigma\ua$ denotes 
the sequence $(a_{n+1}\colon n\in\Z)$.  We use $\overline{\ua}$ to
denote the conjugate sequence $(\overline{a}\colon n\in\Z)$.
In the sequel, $\CA$ denotes a separable subalgebra of $\ell^\infty(\Z)$. 
\subsubsection{The pointed dynamical system associated to an algebra}
\label{subsec:pointeddynamic}
Let  $X$ be 
the {\em Gelfand spectrum} of $\CA$, 
meaning $X$ consists of the set of unitary homomorphisms from $\CA$ 
to the complex numbers.
Letting $\CC(X)$ denote the algebra 
of continuous functions  on $X$, 
we have that there exists an isometric isomorphism of algebras 
$\Phi\colon\CC(X)\to\CA$. For $\ub\in\CA$, the function 
$\Phi\inv(\ub)$
is called \emph{the function associated to} $\ub$.
 
Since $\CA$ is separable, $X$ is a compact metric space. We write 
$d_X$ for a distance on $X$ defining its topology.

The map $\ub\mapsto b_0$ is a character of the algebra $\CA$.   
Thus there exists a point $x_0\in X$ with $f(x_0)=\Phi(f)_0$
for all $f\in\CC(X)$.  
The shift on $\CA$ induces a homeomorphism $T\colon X\to X$ with 
$\Phi(f\circ T)=\Phi(f)\circ\sigma$ for all   $f\in\CC(X)$.
Therefore, for every $f\in\CC(X)$, $\Phi(f)$ is the sequence
$$
 \Phi(f)=\bigl(f(T^nx_0)\colon n\in\Z\bigr)\ .
$$
 In particular, if $f\in\CC(X)$ satisfies $f(T^nx_0) = 0$ for all $n\in\Z$, 
then the sequence given by $\Phi(\ub)=f$ is identically 
zero and so $f$ itself is identically zero.  
It follows that the point $x_0$ 
of $X$ is \emph{transitive}, meaning that its orbit 
$\{ T^nx_0\colon n\in\Z\}$ is dense in $X$.  

We encapsulate this construction in the following definition:
\begin{definition}
The triple $ (X,T,x_0)$ is called the \emph{pointed topological 
dynamical system  associated to the algebra $\CA$}. 
\end{definition}

\subsubsection{Averaging schemes and invariant measures}
\label{subsubsec:averaging}
We first introduce a definition that allows us to average 
any sequence in a subalgebra over a sequence of intervals:
\begin{definition}
Let $\CA$ be a separable subalgebra of $\ell^\infty(\Z)$ and 
$\uI=(I_j\colon j\geq 1)$ be a sequence of intervals whose 
lengths tend to infinity. 
We say that $\uI$ is an \emph{averaging scheme} for $\CA$ if the limit
$$
\LimAvI{\uI}{\ub} :=
\lim_{j\to+\infty}\frac{1}{|I_j|}\sum_{n\in I_j}b_n
$$
exists for all $\ub\in\CA$. 
\end{definition}

Since $\CA$ is separable with respect to the norm 
of $\ell^\infty(\Z)$, for every sequence of intervals whose lengths 
tend to infinity, we can always pass to a subsequence that is an 
averaging scheme for $\CA$. The classical case is when $\uI$ is 
taken to be 
the sequence $([0,j-1]\colon j\geq 1)$, or some subsequence of this 
sequence.

Given an averaging scheme $\uI$ for $\CA$, 
we can associate an invariant probability measure 
$\mu$ on $X$ defined by:
\begin{equation}
\label{eq:Imu}
\int f\,d\mu = \LimAvI{\uI}{f(T^nx_0)}:=
\lim_{j\to+\infty}\frac{1}{|I_j|}\sum_{n\in I_j}f(T^nx_0)
\end{equation}
for all $f\in\CC(X)$.

We claim that all 
ergodic invariant probability measures on $X$ are obtained by this procedure.
Namely, let $\mu$ be such a measure.  
Let $x_1\in X$ be a \emph{generic point} for $\mu$, meaning that 
for all $f\in\CC(X)$, 
$$
\lim_{j\to+\infty}
\frac{1}{j}\sum_{n=0}^{j-1}f(T^nx_1) = \int f\, d\mu \ .
$$
(By the ergodic theorem, $\mu$-almost every point $x_1\in X$ 
is generic.)
Since $x_0$ is a transitive point, there exists a sequence $(k_j\colon j\geq 1)$ 
of integers such that 
$$
\sup_{0\leq n < j}d_X(T^{k_j+n}x_0, T^nx_1) \to 0\text{ as 
}j\to+\infty \ .
$$
So for any continuous function $f$ on $X$, we then have
$$
\lim_{j\to+\infty}
\Bigl(\frac{1}{j}\sum_{n=0}^{j-1}f(T^nx_1) - 
\frac{1}{j}\sum_{n=0}^{j-1}f(T^{k_j+n}x_0)\Bigr) =  0 \ .
$$
Let $\uI$ be the sequence of intervals
$(I_j=[k_j, k_j+j-1]\colon j\geq 1)$. If $\ub\in\CA$ and $f$ is 
the associated function on $X$,  we have
$$
\LimAvI{\uI}{b_n}=\LimAvI{\uI} {f(T^nx_0)}=\int f\,d\mu\ .
$$
Therefore the sequence of intervals $\uI$ 
is an averaging scheme 
for $\CA$ corresponding to the measure $\mu$, and the claim follows.

\subsection{Proofs of properties of the ``seminorms''}
\label{subsec:proofsseminorms}

We use this presentation of the Correspondence Principle 
to derive the properties of the ``seminorms.''  We start 
with the non-negativity that makes the definition possible.
Recall that the bounded sequence $\ua = (a_n\colon n\in\Z)$ satisfies
property $\CP(k)$ on the 
sequence of intervals $\uI$ if for all 
$h=(h_1,\dots,h_k)\in \Z^k$, the limit
$$
 c_h(\uI,\ua)=\LimAvI{\uI}
{\prod_{\epsilon\in\{0,1\}^{k}} C^{|\epsilon|}
a_{n+ \epsilon\cdot h}}
$$
exists.  We show that for a sequence $\ua$ satisfying this, 
the limit
$$
\lim_{H\to+\infty}\frac{1}{H^k}\sum_{h_1, \ldots, h_k=0}^{H-1}
c_h(\uI, \ua)
$$
exists and is non-negative:
\begin{proof}[Proof of Proposition~\ref{prop:seminorms}]
Let $k\geq 2$ be an integer and $\ua=(a_n\colon n\in\Z)$ 
be a bounded sequence that satisfies property $\CP(k)$ on a 
sequence of intervals $\uI$. 
Let $\CA=\CA(\ua)$, $(X,T,x_0)$ be the pointed topological dynamical system 
associated to the algebra $\CA$, and $f\in\CC(X)$ be the function 
associated to the sequence $\ua$. 
Starting with 
the sequence of intervals $\uI$, by passing to a subsequence $\uJ$, 
we extract an averaging scheme for $\CA$. Let $\mu$ be the 
associated measure on $X$. For every $h\in\Z^k$, we have
\begin{equation}
\label{eq:chint}
 c_h(\uI,\ua)=c_h(\uJ,\ua)=\int \prod_{\epsilon\in\{0,1\}^{k}} 
C^{|\epsilon|}f(T^{\epsilon\cdot h}x)\,d\mu(x)\ .
\end{equation}
Let 
$$
 \mu=\int_\Omega \mu_\omega\,dP(\omega)
$$
be the ergodic decomposition of the measure $\mu$. 
The integral~\eqref{eq:chint} can be rewritten as
$$
 \int\Bigl(\int \prod_{\epsilon\in\{0,1\}^{k}} 
C^{|\epsilon|}f(T^{\epsilon\cdot 
h}x)\,d\mu_\omega(x)\Bigr)\,dP(\omega)\ .
$$
By Theorem~\ref{th:cubes}, 
$$
\lim_{H\to+\infty}\frac 1{H^k}\sum_{h_1,\dots,h_k=0}^{H-1}
c_h(\ua,\uI)=\int \nnorm f_{\mu_\omega,k}^{2^k}\,dP(\omega)\ .
$$
Therefore, the announced limit exists and 
is non-negative and we have the statement. 
\end{proof}

Maintaining notation used in the proof, we note that:
\begin{equation}
\label{eq:norminnorm}
 \norm \ua_{\uI,k}= \Bigl(\int\nnorm 
f_{\mu_\omega,k}^{2^k}\,dP(\omega)\Bigr)^{1/2^k}\ .
\end{equation}

We now prove the versions of subadditivity that are satisfied by 
the ``seminorms'':
\begin{proof}[Proof of Propositions~\ref{prop:seminorms2} and~\ref{prop:seminorms3}]
Assume that the bounded 
sequence $\ua$ satisfies properties $\CP(k)$ 
and $\CP(k+1)$ on the sequence of 
intervals $\uI$.  By~\eqref{eq:norminnorm}, the Cauchy-Schwartz 
inequality, and equality~\eqref{eq:k-to-kplus1}, we have
$$
\norm \ua_{\uI,k}^{2^{k+1}} \leq 
\int\nnorm f_{\mu_\omega,k}^{2^{k+1}}\,dP(\omega)
\leq \int\nnorm f_{\mu_\omega,k+1}^{2^{k+1}}\,dP(\omega)
=\norm \ua_{\uI,k+1}^{2^{k+1}}\ .
$$
Thus $\norm \ua_{\uI,k}\leq \norm \ua_{\uI,k+1}$ and 
Proposition~\ref{prop:seminorms3} follows.

Now assume that $\ua$ 
and $\ub$ are bounded sequences and assume that the three sequences $\ua, 
\ub$, and $\ua+\ub$ satisfy property $\CP(k)$ for some sequence of 
intervals $\uI$. We proceed as in the proof of 
Proposition~\ref{prop:seminorms}, taking $\CA$ to be the algebra 
spanned by $\ua$ and $\ub$. If $f$ and $g$ are the functions on $X$ 
associated respectively to $\ua$ and $\ub$, we have that 
\begin{gather*}
 \norm \ua_{\uI,k}^{2^k}  = \int\nnorm 
f_{\mu_\omega,k}^{2^k}\,dP(\omega)\ ;\ 
\norm \ub_{\uI,k}^{2^k}  =\int \nnorm 
g_{\mu_\omega,k}^{2^k}\,dP(\omega)\ ; \\
\norm {\ua+\ub}_{\uI,k}^{2^k}  = \int \nnorm 
{f+g}_{\mu_\omega,k}^{2^k}\,dP(\omega)\ .
\end{gather*}
Therefore 
$$
 \norm{\ua+\ub}_{\uI,k}\leq\norm\ua_{\uI,k}+\norm\ub_{\uI,k}
$$
and Proposition~\ref{prop:seminorms2} follows. 
\end{proof}

\subsection{A Cauchy-Schwartz-Gowers type result}
We have an inequality similar to that satisfied by the Gowers norms in the
finite setting and by the HK-seminorms, as given in~\eqref{eq:HKCSG}:
\begin{proposition}
\label{prop:HKCSG}
For every $\epsilon\in\{0,1\}^k$, let 
$\ua(\epsilon)=\bigl(a_n(\epsilon)\colon n\in\Z\bigr)$ be a bounded 
sequence. 
Let $\uI$ be  a sequence of intervals whose lengths tend to 
infinity such that
$$
c_h:= \lim_{j\to+\infty}\frac 1{|I_j|}\sum_{n\in I_j}
{\prod_{\epsilon\in\{0,1\}^k} a_{n+\epsilon\cdot h}}
$$
exists for every $h\in\Z^k$. Then the limit
$$
\lim_{H\to+\infty} \frac 1{H^k}\sum_{h_1,\dots,h_k=0}^{H-1}c_h
$$
exists.

Moreover, if all the sequences $\ua(\epsilon)$ satisfy property 
$\CP(k)$ on $\uI$, then 
\begin{equation}
\label{ineq:HKCSG}
\Bigl|\lim_{H\to+\infty} \frac 
1{H^k}\sum_{h_1,\dots,h_k=0}^{H-1}c_h\Bigr|
\leq\prod_{\epsilon\in\{0,1\}^k}
\norm{\ua(\epsilon)}_{\uI,k}\ .
\end{equation}
\end{proposition}

\begin{proof}
The proof of the convergence is similar to the proof of 
Proposition~\ref{prop:seminorms}, but we set $\CA$ to be the 
algebra spanned by the $2^k$ sequences $\ua(\epsilon)$, 
$\epsilon\in\{0,1\}^k$. Maintaining notation as that proof, 
for every $\epsilon\in\{0,1\}^k$ we let $f_\epsilon$ denote the function 
associated to the sequence $\ua(\epsilon)$.  
It follows from inequality~\eqref{eq:HKCSG} that 
\begin{multline*}
\Bigl|\lim_{H\to+\infty} \frac 
1{H^k}\sum_{h_1,\dots,h_k=0}^{H-1}c_h\Bigr|=\Bigl|\int\Bigl(\int 
\prod_{\epsilon\in\{0,1\}^k } 
f_\epsilon (x_\epsilon)\,
d\mu_\omega\type k(\ulx)\Bigr)\,dP(\omega)\Bigr| \\
\leq\int \prod_{\epsilon\in\{0,1\}^k } \nnorm{f_\epsilon}_{\mu_\omega,k}
\,dP(\omega)
\leq\prod_{\epsilon\in\{0,1\}^k }\Bigl(
\int \nnorm{f_\epsilon}_{\mu_\omega,k}^{2^k}\,dP(\omega)\Bigr)^{1/2^k}
= \prod_{\epsilon\in\{0,1\}^k} \norm{\ua(\epsilon)}_{\uI,k}\ .
\end{multline*}
\end{proof}

Using 
relations~\eqref{eq:k-to-kplus1} and~\eqref{eq:norminnorm}, we deduce that:
\begin{proposition}
\label{prop:kk+1sequences}
Assume that the bounded sequence $\ua$ satisfies property 
$\CP(k+1)$ on $\uI$.
Then
$$
\lim_{H\to+\infty}\frac 1H\sum_{h=0}^{H-1}
\norm{\sigma^h\ua.\overline{\ua}}_{\uI,k}^{2^k}=
\norm\ua_{\uI,k+1}^{2^{k+1}}\ .
$$
\end{proposition}
Note that the hypothesis implies that for every integer $h\geq 1$, the 
sequence $\sigma^h\ua.\overline{\ua}$ satisfies property $\CP(k)$ 
on $\uI$.

\subsection{The uniformity seminorms}
We also use the Correspondence Principle to derive properties of the 
uniformity seminorms:
\begin{proposition}
\label{prop:supergodic}
Let $k\geq 1$ be an integer, $\ua$ be a bounded sequence, 
$(X,T,x_0)$ the associated 
pointed dynamical system, and $f\in\CC(X)$ be the function associated to 
$\ua$.  Then
$$
\norm\ua_{U(k)}=
\sup_{\mu\text{ \upshape ergodic}}\nnorm f_{\mu,k}\ ,
$$
where the supremum is taken over all ergodic measures $\mu$ on $X$.
\end{proposition}

\begin{proof}
It follows from~\eqref{eq:norminnorm} 
that if we raise the left hand side to the power $2^k$, then 
it is bounded by the right hand side raised to the power $2^k$.  
Conversely, in Section~\ref{subsubsec:averaging} we showed that
every ergodic measure $\mu$ on $X$ is associated to an averaging 
scheme $\uI$ for the algebra $\CA(\ua)$.  By applying~\eqref{eq:norminnorm} 
again, we have that 
$\nnorm f_{\mu,k}=\norm\ua_{\uI,k}\leq\norm\ua_{U(k)}$.
\end{proof}

Proposition~\ref{prop:uk-seminorm} follows immediately; it could also 
be derived directly from Proposition~\ref{prop:seminorms2}.

\begin{remark}
We note that there are important differences between the 
uniformity seminorms and the HK-seminorms. For 
example,
the formula given by Proposition~\ref{prop:kk+1sequences} comes from, 
and is similar to, formula~\eqref{eq:k-to-kplus1} for the 
HK-seminorms. 
We deduce that 
$$
\norm \ua_{U(k+1)}^{2^{k+1}}
\leq\liminf_{H\to+\infty}
\frac 1H\sum_{h=0}^{H-1}
\norm{\bar \ua.\sigma^h\ua}_{U(k)}^{2^k}\ .
$$
But in general,
 the $\liminf$ on the right hand side of this 
equation is not a
limit and equality does not hold.
\end{remark}

\section{A duality in nilmanifolds and direct results}
\label{sec:duality}
\subsection{Measures and norms for nilsystems}

Throughout 
 this section, we assume 
that $k\geq 2$ 
is an integer and $(X=G/\Gamma,\mu,T)$ is an ergodic $(k-1)$-step nilsystem, where 
$T$ is the translation by $\tau\in G$.  
As explained in Section~\ref{sec:tools}, 
we reduce to  the case that $G$ is spanned by its 
connected component $G_0$ of the identity and by $\tau$.

We review properties of the measure $\mu\type k$ and of 
the seminorm $\nnorm\cdot_k$ in 
this particular case.  
Most of these properties are established in~\cite{HK} 
or~\cite{GT}, 
but often in a very different context and with very 
different terminology from that used here.  We include some
proofs for completeness, but as they are far from the main 
topics of the article, we defer them to 
Appendix~\ref{appendix:cubes}. 
This appendix also includes some 
properties we need that are not stated elsewhere. 

We use the notation for $2^k$-Cartesian powers introduced in 
Section~\ref{sec:tools}.  We summarize the properties that we need:

\begin{theorem}
\label{th:summary}

\strut
\begin{enumerate}
\item
\label{it:Haarnil}
The measure $\mu\type k$ is the Haar measure of a sub-nilmanifold $X_k$ 
of $X\type k$. The transformations $T\type k$ and $T\type k_i$, 
$1\leq i\leq k$, act on $X_k$ by translation 
and  $X_k$ is 
 ergodic (and thus uniquely ergodic and minimal) under these transformations.
\item
\label{it:Phi}
Let $X_{k*}$ be the image of $X_k$ under the 
projection $\ulx\mapsto\ulx_*$ from $X\type k$ to $X\type k_*=X^{2^k-1}$.
There exists a smooth map $\Phi\colon X_{k*}\to X_k$ such that 
$$
 X_k=\bigl\{(\Phi(\ulx_*),\ulx_*)\colon \ulx\in X_{k*}\bigr\}\ .
$$
\item
\label{it:norm}
$\nnorm\cdot_k$ is a norm on $\CC(X)$.
\item
\label{it:Wx}
For every $x\in X$, let $W_{k,x}=\{ \ulx\in X_k\colon x_\bzero=x\}$.
Then $W_{k,x}$ is uniquely ergodic under the transformations 
$T\type k_i$, $1\leq i\leq k$.
\item\label{it:measureWx}
For every $x\in X$, let $\rho_x$ be the invariant measure of $W_{k,x}$. 
Then for every $x\in X$ and $g\in G$, $\rho_{g.x}$ is the image of 
$\rho_x$ under the translation by $g\type k=(g,g,\dots,g)$.
\end{enumerate}
\end{theorem}

The nilmanifold $X_k$ is 
defined independently of the transformation $T$ and it 
only depends on the structure of the nilmanifold $X$. 
This implies that the measure $\mu\type k$ and the norm 
$\nnorm\cdot_k$ do not depend on the 
transformation $T$ on $X$, provided that $T$ is an ergodic transformation.
These are geometric, and not dynamical, objects.

\subsection{Uniform convergence}
Using part~\eqref{it:Wx} of Theorem~\ref{th:summary} we deduce:
\begin{corollary}
\label{cor:uniform}
Let $f_\epsilon$, $\epsilon\in\{0,1\}^k_*$ 
be $2^k-1$ continuous functions on $X$. 
For every $x\in X$ we 
have 
$$
\frac 1{H^k}\sum _{h_1,\dots,h_k=0}^{H-1}
\prod_{\epsilon\in\{0,1\}^k_*}f_\epsilon(T^{\epsilon\cdot h}x)
\to  \int 
\prod_{\epsilon\in\{0,1\}^k_*}f_\epsilon(x_\epsilon)\,d\rho_x(\ulx)
$$
as $H\to+\infty$.  
 Moreover, the convergence is uniform in $x\in X$.
\end{corollary}

\begin{proof}
The corollary  follows easily from part~\eqref{it:Wx} of 
Theorem~\ref{th:summary} by a classical 
argument. 
Let $(x_j\colon j\geq 1)$ be a sequence in $X$ converging to some $x\in 
X$ and let $(H_j\colon j\geq 1)$ be a sequence of integers tending 
to infinity.

For every $j$, let $\nu_j$ be the measure
$$
\nu_j:= \frac 1{H_j^k}\sum _{h_1,\dots,h_k=0}^{H_j-1}
\underset{\epsilon\in\{0,1\}^k}{\otimes}\delta_{T^{\epsilon\cdot h}x_j} 
$$
on $X\type k$ and let  $\nu$ be any weak limit of this sequence of measures. For 
every $j$, the measure $\nu_j$ is concentrated on $W_{k,x_j}$. 
Since $X_k$ 
is closed in $X\type k$, the measure $\nu$ is concentrated on $W_{k,x}$.
Moreover, for every $j$ and for $1\leq i\leq k$, the difference 
between the measures $\nu_j$ and $T\type 
k_i\nu_j$ are at a distance 
$\leq 2/H_j$ in the norm of total variation.  It follows that 
$ \nu$ is invariant under $T\type k_i$ for $i=1, \ldots, k$.  
By unique ergodicity of $W_{k,x}$, 
we have that $\nu$ is the invariant measure $\rho_x$ of $W_{k,x}$. 

We  have shown that 
the sequence $(\nu_j\colon j\geq 1)$ 
of measures converges weakly to the measure 
$\rho_x$.
It follows that if $f_\epsilon$, $\epsilon\in\{0,1\}^k_*$, are 
continuous functions on $X$, then
\begin{align*}
\frac 1{H_j^k}\sum _{h_1,\dots,h_k=0}^{H_j-1}
\prod_{\epsilon\in\{0,1\}^k_*} f_\epsilon(T^{\epsilon\cdot h}x_j)
& =\int \prod_{\epsilon\in\{0,1\}^k_*} f_\epsilon(x_\epsilon)\,d\nu_j(\ulx)\\
& \to 
\int 
\prod_{\epsilon\in\{0,1\}^k_*}f_\epsilon(x_\epsilon)\,d\rho_x(\ulx)
\end{align*}
as $j\to+\infty$ and the result follows.
\end{proof}

We apply this result when $f$ is a continuous function on 
$X$ and $f_\epsilon=C^{|\epsilon|}f$ for every 
$\epsilon\in\{0,1\}^k_*$.  From Corollary~\ref{cor:dual}, we have 
that the averages in 
Corollary~\ref{cor:uniform} converge in $L^2(\mu)$ to the function 
$\CD_kf$.  Therefore:

\begin{corollary}
\label{cor:uniform2}
Let  $f$ be a continuous function on $X$. Then
$$
 \CD_kf(x)=\int \prod_{\epsilon\in\{0,1\}^k_*}C^{|\epsilon|}f(x_\epsilon)\,d\rho_x(\ulx)
$$
and the function $\CD_kf$ is the uniform limit of the sequence
$$
 \frac 1{H^k}\sum _{h_1,\dots,h_k=0}^{H-1}
\prod_{\epsilon\in\{0,1\}^k_*}f_\epsilon(T^{\epsilon\cdot h}x)\ .
$$
Thus $\CD_kf$ is a continuous function on $X$.
\end{corollary}

In particular, the function $\CD_kf$ is a geometric object: it 
does not depend on the transformation $T$ on $X$.

\begin{corollary}
\label{cor:smooth}
If $f$ is a smooth function on $X$, then $\CD_k f$ is a smooth 
function on $X$.
\end{corollary}
\begin{proof}Let $x_0\in X$. Then, by Corollary~\ref{cor:uniform2} 
and part~\eqref{it:measureWx} of Theorem~\ref{th:summary}, for every 
$g\in G$ we have
$$
 \CD_kf(g.x_0)=\int 
\prod_{ \epsilon\in\{0,1\}^k_* }C^{|\epsilon|}f(g.x_\epsilon)\,d\rho_{x_0}(\ulx)\ .
$$
Thus the function $g\mapsto\CD_kf(g.x_0)$ is a smooth function on 
$G$ and the result follows.
\end{proof}

\begin{remark}
Let $x\in X$. 
Since the measure $\rho_x$ is invariant under the 
transformations $T\type k_i$, it follows that the image of this 
measure under the projection $\ulx\mapsto x_\epsilon$  
for every $\epsilon\in\{0,1\}^k$ is invariant 
under $T$ and thus is equal to $\mu$. Therefore if $f_\epsilon$, 
$\epsilon\in\{0,1\}^k_*$, are  continuous 
functions on $X$, the H\"older inequality gives:
$$
\Bigl|
\int\prod_{\epsilon\in\{0,1\}^k_*}f_\epsilon(x_\epsilon)\,d\rho_x(\ulx)\Bigr|
\leq\prod_{\epsilon\in\{0,1\}^k_*} 
\norm{f_\epsilon}_{L^{2^k-1}(\mu)}\ .
$$
By density we deduce that for every $f\in L^{2^k-1}(\mu)$ the 
function $\CD_kf$ is continuous on $X$ and that 
$$\norm{\CD_kf}_\infty\leq\norm f_{L^{2^k-1}(\mu)}^{2^k-1}\ .$$
\end{remark}

\subsection{The dual norm}
\label{subsec:dualnorm}
\begin{definition}
Let the space $\CC(X)$ of continuous functions on $X$ be endowed with 
the norm $\nnorm\cdot_k$. 
Since $\nnorm f_k\leq\norm f_{L^{2^k}(\mu)}$ for every $f\in\CC(X)$,
the dual of this space can be identified with a subspace of 
$L^{2^k/(2^k-1)}(\mu)$. 
We call this space {\em the dual space} and denote it by $\CC(X)_k^*$.  
We write  $\nnorm h_k^*$ for the {\em dual norm} of a function 
$h\in\CC(X)^*_k$.
\end{definition}
In other words, a function 
$h\in L^{2^k/(2^k-1)}(\mu)$ belongs 
to the dual space $\CC(X)^*_k$ if there exists a constant $C$ with
\begin{equation}
\label{eq:dualnorm}
\Bigl|\int f\,h\,d\mu\Bigr|\leq C\,\nnorm f_k
\end{equation}
 for every $f\in \CC(X)$
and $\nnorm f_k^*$ is the smallest constant $C$ with this property.

We note that the  dual space and the dual norm  $\nnorm 
\cdot _k^*$  are geometric, not dynamical, objects.

We give two methods to build functions in the dual space.
Let $f$ be a function on $X$, belonging to $L^{2^k}(\mu)$. By 
characterization~\eqref{eq:defDk}  of the dual function and
inequality~\eqref{eq:HKCSG}, we have that 
for every $h\in\CC(X)$, 
$$
\Bigl|\int h.\CD_k f\,d\mu\Bigr|\leq \nnorm h_k\,\nnorm f_k^{2^k-1}\ .
$$ 
Thus $\CD_k f$ belongs to the dual space  and 
$\nnorm{\CD_k f}_k^*\leq \norm f_k^{2^k-1}$. On the other hand,
$$
 \nnorm f_k\,\nnorm{\CD_k f}_k^*\geq \int  f\,\CD_k f\,d\mu=\nnorm 
f_k^{2^k}
$$
and we conclude that
\begin{equation}
\label{eq:Dkfdualfunction}
\nnorm{\CD_k f}_k^*=\nnorm f_k^{2^k-1}\ .
\end{equation}

We now show:
\begin{proposition}
\label{prop:dualsmooth}
The dual space  $\CC(X)^*_k$ contains all 
smooth functions on $X$. 
\end{proposition}

\begin{proof}
Let $f$ be a smooth function on $X$ and let $X_{k*}$ and $\Phi$ be 
the set and the map defined 
in part~\eqref{it:Phi} of Theorem~\ref{th:summary}.

Then $f\circ\Phi$ is a smooth function on  $X_{k*}$ and there 
exists a smooth function $F$ on $X\type k_*$ whose restriction to 
$X_{k*}$ is equal to $f\circ\Phi$.
This function can be written as
$$
 F(\ulx_*)=\sum_{j=1}^\infty\prod_{\epsilon\in\{0,1\}^k_*}
f_{j,\epsilon}(x_\epsilon)\ ,
$$
where  the functions $f_{j,\epsilon}$, $j\geq 1$ and
$\epsilon\in\{0,1\}^k_*$, are continuous functions on $X$ satisfying
$$
\sum_{j= 1}^\infty
\prod_{\epsilon\in\{0,1\}^k_*}\norm{f_{j,\epsilon}}_\infty 
<+\infty\ .
$$
For every continuous function $h$ on $X$, we have 
\begin{align*}
\Bigl|\int f\,h\,d\mu\Bigr|
= &\Bigl|\int h(x_\bzero)\,.\,f\circ\Phi(\ulx_*)\,d\mu\type k(\ulx)\Bigr|
=  \Bigl|\int h(x_\bzero)\,.\,F(\ulx_*)\,d\mu\type k(\ulx)\Bigr|\\
\leq & \sum_{j=1}^\infty \Bigl|\int h(x_\bzero)\,\prod_{\epsilon\in\{0,1\}^k_*}
h_{j,\epsilon}(x_\epsilon)\,d\mu\type k(\ulx)\Bigr|\\
\leq & \sum_{j=1}^\infty \nnorm 
h_k\prod_{\epsilon\in\{0,1\}^k_*}\nnorm{h_{j,\epsilon}}_k\\
\leq &
\nnorm  h_k 
\sum_{j=1}^\infty \prod_{\epsilon\in\{0,1\}^k_*}\norm{h_{j,\epsilon}}_\infty \ .
\end{align*}
where the next to last inequality follows from~\eqref{eq:HKCSG}.  
The announced statement follows.  
\end{proof}

A similar proof is used in~\cite{GT} in the 
finite setting.

The hypothesis of smoothness is too strong and could be replaced by 
weaker assumptions. It is probably sufficient to assume that $f$ is 
Lipschitz with respect to some smooth metric on $X$.
Computing the dual norm of $f$, or even bounding it in an explicit
way seems to be difficult. The regularity of the map $\Phi$ 
should play a role, but in order to  define this, we would first 
need to choose a metric on $X$.

\begin{proposition}
\label{prop:dualspan}
The unit ball of $\CC(X)_k^*$ is 
the closure in $L^{2^k/(2^k-1)}(\mu)$ of the convex hull of the set
$$
 \{ \CD_k f\colon f\in\CC(X),\ \nnorm f_k\leq 1\}\ .
$$
\end{proposition}

\begin{proof}
Let $B$ be the set in the statement. 
By~\eqref{eq:Dkfdualfunction}, for $f\in\CC(X)$ with $\nnorm f_k\leq 1$, we have 
that $\CD_k f$ belongs to the unit ball of $\CC(X)^*_k$. 
Since this ball is closed in the  norm of $L^{2^k/(2^k-1)}(\mu)$, it 
contains $B$.

On the other hand, let $f$ be a nonzero function belonging to $L^{2^k}(\mu)$
 and let $h=\nnorm f_k\inv. f$.  As the map $\CD_k\colon 
L^{2^k}(\mu)\to L^{2^k/(2^k-1)}(\mu)$ is continuous, by density we 
have that $\CD_k h\in B$. As $\int 
f.\CD_k h\,d\mu=\nnorm f_k$, the Hahn-Banach Theorem gives the opposite 
inclusion.
\end{proof}

\subsection{Direct theorem (upper bound) }
\label{subsec:ub}
We now have assembled the ingredients to prove 
Theorem~\ref{th:upperbound}. As we have
not yet defined the norm $\nnorm\ub_k^*$ of a smooth nilsequence 
$\ub$, we state this theorem in a modified version. 
\begin{theorem*}[Modified Direct Theorem]
Let $\ua=(a_n\colon n\in\Z)$ be a bounded sequence that satisfies property $\CP(k)$ on the 
sequence of intervals $\uI = (I_j\colon j\geq 1)$.  Let $(X,T,\mu)$ 
be an ergodic $(k-1)$-step nilsystem, $x_0\in X$, and $f$ be a smooth 
function on $X$. 
Then 
$$
\limsup_{j\to+\infty}\Bigl|
\frac 1{|I_j|}\sum_{n\in I_j}a_n\,f(T^nx_0)\Bigr|
\leq\norm\ua_{\uI,k}\,\nnorm f_k^*\ .
$$
\end{theorem*}
 
\begin{proof}
\strut
\subsubsection{}
We begin with the case that 
$ f=\CD_k \phi$ for some continuous function $\phi$ on $X$ with 
$\nnorm\phi_k=1$.

By substituting a subsequence for $\uI$, we can assume that 
for every $h=(h_1,\dots,h_k)\in\Z^k$,  the averages on $I_j$ of
$$
a_n \prod_{\epsilon\in\{0,1\}^k_*}C^{|\epsilon|}\phi(T^{n+\epsilon\cdot 
h}x_0)
$$
converge.

Fix $\delta>0$. By Corollary~\ref{cor:uniform2}, for every 
sufficiently large $H$ we have that 
$$
 \Bigl|\frac 1{H^k}\sum_{h_1,\dots,h_k=0}^{H-1} 
a_n\prod_{\epsilon\in\{0,1\}^k_*} 
C^{|\epsilon|}\phi(T^{n+\epsilon\cdot h}x_0)-a_nf(T^nx_0)\Bigr|<\delta
$$
for every $n\in\Z$ and so
$$
\Bigl|\frac 1{H^k}\sum_{h_1,\dots,h_k=0}^{H-1} \Bigl(\frac 
1{|I_j|}\sum_{n\in I_j} a_n\prod_{\epsilon\in\{0,1\}^k_*} 
C^{|\epsilon|}\phi(T^{n+\epsilon\cdot h}x_0)\Bigr)- \frac 
1{|I_j|}\sum_{n\in I_j} a_nf(T^nx_0)\Bigr|<\delta\ .
$$
for every $j\geq 1$. Taking the limit as $j\to+\infty$ along a subsequence, 
we have that for every sufficiently large $H$,
\begin{multline*}
 \LimsupAvI{\uI}{a_nf(T^nx_0)}\\
\leq\delta+ 
\left|\frac 1{H^k}\sum_{h_1,\dots,h_k=0}^{H-1}\LimAvI{\uI}
{a_n\prod_{\epsilon\in\{0,1\}^k_*} 
C^{|\epsilon|}\phi(T^{n+\epsilon\cdot h}x_0)}\right|\ .
\end{multline*}
We conclude that 
\begin{multline*}
 \LimsupAvI{\uI}{a_nf(T^nx_0)}\\
\leq
\left|\lim_{H\to+\infty} 
\frac 1{H^k}\sum_{h_1,\dots,h_k=0}^{H-1}\LimAvI{\uI}
{a_n\prod_{\epsilon\in\{0,1\}^k_*} 
C^{|\epsilon|}\phi(T^{n+\epsilon\cdot h}x_0)}\right|\ .
\end{multline*}
The existence of the limit for $H\to+\infty$  is given by 
Proposition~\ref{prop:HKCSG}. Using 
Inequality~\eqref{ineq:HKCSG}
 and Corollary~\ref{cor:normtotergo}, we 
have that the last quantity is bounded by
$$
\norm\ua_{\uI,k}\,.\, \norm{ (\phi(T^nx_0)\colon n\in \Z)}_{\uI,k}^{2^k-1}
=\norm\ua_{\uI,k}\,.\, \nnorm\phi_k^{(2^k-1)/2^k}=  \norm\ua_{\uI,k}\ .
$$

\subsubsection{} We now turn to the general case. We can assume that 
$\nnorm f_k^*\leq 1$.

Fix $\delta >0$. 
By Proposition~\ref{prop:dualsmooth},
we can write $f=f_1+f_2$, where $f_1$ is a convex combination of 
functions considered in the first part and 
$\norm{f_2}_{L^{2^k/(2^k-1)}(\mu)}<\delta$.  
The contribution of $f_1$ to the $\limsup$ of the 
averages is bounded by $1$. 

For every $j\geq 1$, by the 
H\"older inequality we have
$$
\Bigl|\frac 1{|I_j|}\sum_{n\in I_j}a_nf_2(T^nx_0)\Bigr|\leq
\norm\ua_\infty 
\Bigl(\frac 1{|I_j|}\sum_{n\in 
I_j}|f_2(T^nx_0)|^{2^k/(2^k-1)}\Bigr)^{(2^k-1)/2^k}\ .
$$
Since both $f$ and $f_1$ are 
continuous, so is $f_2$. Therefore, by unique 
ergodicity of $(X,T)$, the averages of $|f_2(T^nx_0)|^{2^n/(2^n-1)}$ 
converge to the integral of the function $|f|^{2^n/(2^n-1)}$ and we have 
that 
$$\LimsupAvI{\uI}{a_nf_2(T^nx_0)}\leq\delta\ .$$
The result follows.
\end{proof}

\subsection{The dual norm for smooth nilsequences}
\begin{corollary}
\label{cor:normsequence}
Let $(X,\mu,T)$ and $(Y,\nu,S)$ be ergodic $(k-1)$-step nilsystems, 
$x_0\in X$, $y_0\in Y$, $f$ be a smooth function on $X$, and $g$ a smooth 
function on $Y$. If $f(T^nx_0)=g(S^ny_0)$ for every $n\in\Z$, then 
$\nnorm f_{\mu,k}^*=\nnorm g_{\nu,k}^*$.
\end{corollary}

\begin{proof}
Fix $\delta>0$. By definition of $\nnorm f_{\mu,k}^*$, there exists a 
continuous function $h$ on $X$ with 
$$
 \nnorm h_{\mu,k}=1\text{ and }\Bigl|\int f\,h\,d\mu\Bigr|\geq
\nnorm f_{\mu,k}^*-\delta\ .
$$
By unique ergodicity of $X$,
$$
\Bigl|\int f\,h\,d\mu\Bigr|
=\lim_{N\to+\infty}\Bigl|\sum_{n=0}^{N-1} 
f(T^nx_0)h(T^nx_0) \Bigr|=
\lim_{N\to+\infty}\Bigl|\sum_{n=0}^{N-1} 
g(S^ny_0)h(T^nx_0) \Bigr| \ .
$$
Let $\uI$ be the sequence of intervals $(I_N = [0,N-1]\colon N\geq 1)$.
By Corollary~\ref{cor:normtotergo}, the sequence $(h(T^nx_0)\colon 
n\in\Z)$ satisfies property $\CP(k)$ on $\uI$ and 
$\norm{(h(T^nx_0)\colon n\in\Z)}_{\uI,k}=\nnorm h_{\mu,k}=1$.  By
the Modified Direct Theorem, we have that
$$
 \lim_{N\to+\infty}\Bigl|\sum_{n=0}^{N-1} 
g(S^ny_0)h(T^nx_0) \Bigr|\leq\nnorm g_{\nu,k}^*
$$
and so 
$\nnorm f_{\mu,k}^*-\delta \leq \nnorm g_{\nu,k}^*$.  
Exchanging the roles of $f$ and $g$, we obtain the announced equality.
\end{proof}

Using this corollary, we define:
\begin{definition}
Let $\ub$ be a $(k-1)$-step smooth nilsequence. We define
$\nnorm\ub_k^*=\nnorm f_{\mu,k}^*$, where $f$ is a smooth function 
on an ergodic $(k-1)$-step nilsystem $(X,\mu,T)$
and $x_0\in X$ is chosen such that  $b_n=f(T^nx_0)$ for every $n$.
\end{definition}
Using this definition, the Direct Theorem (Theorem~\ref{th:upperbound}) 
is a reformulation of  the Modified Direct Theorem 
of Section~\ref{subsec:ub}.  

\subsection{The case $k=2$}
Let $X$ be a $1$-step nilmanifold, that is, a compact abelian Lie 
group, and let $f$ be a smooth function on $X$. Let $\wh X$ be the dual group 
of $G$. Then the Fourier series of $f$ is
$$
 f(x)=\sum_{\chi\in\wh X}\wh f(\chi)\,\chi(x)\ , \text{ where 
}\sum_{\chi\in\wh X}|\wh f(\chi)|<+\infty\ .
$$
An easy computation using the definition gives
$$
 \nnorm f_2=\Bigl(\sum_{\chi\in\wh X}|\wh f(\chi)|^4\Bigr)^{1/4}\ .
$$
Therefore we have 
$$
\nnorm f_2^*=\Bigl(\sum_{\chi\in\wh X}|\wh f(\chi)|^{4/3}\Bigr)^{3/4}\ .
$$
If $T$ is an ergodic translation on $X$, $x_0\in X$, and $\ub$ is the 
sequence given by
$b_n=f(T^nx_0)$ for every $n$, we recover the formula for 
$\nnorm\ub_2^*$ given in Section~\ref{subsec:kequals2} and 
Proposition~\ref{prop:upperkequals2}.

\subsection{Some convergence results}

\begin{corollary}
\label{cor:decomp-converges}
Let $k\geq 2$ be an integer, $\uI=(I_j\colon j\geq 1)$ 
be a sequence of intervals whose 
lengths tend to infinity, and 
let $\ua=(a_n\colon n\in\Z)$ be a bounded sequence.  
Assume that for every $\delta>0$, 
there exists a $(k-1)$-step nilsequence $\ua'$ such that 
the sequence $\ua-\ua'$ satisfies property $\CP(k)$ on $\uI$ and 
$\norm{\ua-\ua'}_{\uI,k}<\delta$.
Then for every $(k-1)$-step nilsequence $\ub=(b_n\colon n\in\Z)$, 
the limit
$$
\lim_{j\to+\infty}\frac{1}{|I_j|}\sum_{n\in I_j}a_nb_n
$$
exists.
\end{corollary}
\begin{proof}
By density, we can restrict to the case that $\ub$ is a smooth 
nilsequence. Let $\delta>0$ and 
the nilsequence $\ua'$ be as in the statement. Since the product 
sequence $\ua'\ub$ is a nilsequence, its
averages converge. By Theorem~\ref{th:upperbound},
$$
\LimsupAvI{\uI}{(a_n-a'_n)b_n}\leq\delta\nnorm b_k^* \ .
$$
It follows that the averages on $I_j$ of $a_nb_n$ form a Cauchy 
sequence.
\end{proof}
By the same argument, we have:
\begin{corollary}
\label{cor:decomp-converges2}
Let $k\geq 2$ be an integer and 
$\ua=(a_n\colon n\in\Z)$ 
be a bounded sequence. Assume that for every $\delta>0$, 
there exists a $(k-1)$-step nilsequence $\ua'$ such that 
$\norm{\ua-\ua'}_{U(k)}<\delta$.
Then for every $(k-1)$-step nilsequence $\ub=(b_n\colon n\in\Z)$, 
the averages of the sequence $a_nb_n$ converge, meaning 
that
the limit
$$
\lim_{j\to+\infty}\frac{1}{|I_j|}\sum_{n\in I_j}a_nb_n
$$
exists for all sequences of intervals $\uI = (I_j\colon j\geq 1)$ 
whose lengths tend to infinity. 
\end{corollary}
This Corollary is the direct implication of 
Theorem~\ref{th:cond-conv}. Propositions~\ref{prop:genpol} 
and~\ref{prop:morse} provide 
examples of sequences satisfying the hypothesis of this Corollary.

\section{The correspondence principle revisited and  inverse 
theorems}
\label{sec:extension}

\subsection{An extension of the  correspondence principle}

We recall that a topological dynamical system $(Y,S)$ is \emph{distal} if 
for every $y,y'\in Y$ with $y\neq y'$, then 
$$
\inf_{n\in\Z}d_Y(T^ny,T^ny')>0
$$
where $d_Y$ denotes a distance defining the topology of $Y$.
\begin{proposition}
\label{prop:extended}
Let $(X,T)$ be a topological dynamical system, $x_0\in X$ a 
transitive point, and $\mu$ an invariant ergodic measure on $X$. 
Let $(Y,S)$ be a distal topological dynamical system, $\nu$ an 
invariant measure on $Y$, and $\pi\colon (X,\mu,T)\to(Y,\nu,S)$ a 
measure theoretic factor map.

Then there exist a point $y_0\in Y$ and a sequence of intervals 
$\uI=(I_j\colon j\geq 1)$ whose lengths tend to infinity such that 
for every continuous function 
$f$ on $X$ and every continuous function $g$ on $Y$, 
$$
 \int f(x).g\circ\pi(x)\,d\mu(x)=
\lim_{j\to+\infty}\frac 1{|I_j|}\sum_{n\in I_j}
f(T^nx_0).g(S^ny_0)\ .
$$
\end{proposition}

If the system $(X,T)$ and the point $x_0$ are associated to a 
sequence as in  Section~\ref{sec:correspondence} and if $Y$ 
denotes the Kronecker factor of 
$(X,\mu,T)$, then the sequence of intervals $\uI$ given by the 
Proposition plays the same role as the  ``Kronecker 
complete processes'' of~\cite{BFW}. Our construction is (we hope) simpler
and  works in a more 
general setting: below we use it when $Y$ is a nilsystem.

 \begin{proof}
We write $d_X(\cdot,\cdot)$ and $d_Y(\cdot,\cdot)$ for distances on 
$X$ and $Y$ defining the topologies of these spaces.

\subsubsection{Construction of an extension of $X$}
Let $\CB$ be the closed (in norm) subalgebra of $L^\infty(\mu)$ 
that is spanned by $\CC(X)$ and the functions $g\circ\pi$ with 
$g\in\CC(Y)$.  
This algebra is unitary, separable, and invariant under complex conjugation 
and under $T$.  

Let $W$ be the Gelfand spectrum of this algebra.
Since $\CB$ is separable, $W$ is a compact metrizable space.  
  By definition, there exists  an isometric isomorphism of algebras $\Psi\colon \CC(W)\to\CB$.

As in Section~\ref{subsec:pointeddynamic},
there exists a homeomorphism $R\colon W\to W$ 
satisfying $\Psi(f\circ T) = \Psi(f)\circ R$ for all functions $f\in\CC(W)$.  

The inclusion of $\CC(X)$ in $\CB$ induces a continuous surjective map 
$p\colon W \to X$ satisfying $f\circ p = \Psi(f)$ for every continuous function $f$ on $X$ and we have that 
$T\circ p = p \circ R$.  
Similarly, the map $g\mapsto g\circ\pi$ from $\CC(Y)$ to 
$\CB$ is an isometric homomorphism of algebras and thus induces a continuous surjective map $q\colon W\to Y$ 
satisfying $g\circ q = \Psi(g\circ\pi)$ for all continuous functions $g$ on $Y$.  
We have that $S\circ q = q \circ R$.  So, $p\colon (W,R)\to(X,T)$ and 
$q\colon(W,R)\to(Y,S)$ are factor maps, in the topological sense.

The map $f\mapsto \int f\,d\mu$ is a positive linear form on the 
algebra $\CB$ and 
thus there exists a unique probability measure $\rho$
on $W$  satisfying 
$$
\int f\, d\mu = \int\Psi(f)\,d\rho \text{ for all functions } f\in\CB \ .
$$

 Since 
$\Psi(f\circ T) = \Psi(f)\circ R$ for all $f\in\CB$ and $\mu$ is 
invariant under $T$, 
the measure $\rho $ is invariant under $R$.
Since $\Psi(f) = f\circ p$ for all continuous functions $f$ on $X$, 
we have that
the image of $\rho $ under $p$ is equal to $\mu$. Therefore, $p\colon 
(W,\rho,R)\to(X,\mu,T)$ is a measure theoretic factor map.
  Moreover, for every function 
$f\in\CB$, 
$$
\int |\Psi(f)|^2\,d\rho = 
\int \Psi(|f|^2)\,d\rho = 
\int |f|^2\,d\mu 
$$ 
and the map $\Psi$ is an isometry from the space $\CB$ endowed with the 
norm $L^2(\mu)$ into the space $L^2(\rho)$.  
Since $\CC(X)$ is dense in $\CB$ under the $L^2(\mu)$ norm and 
since $\Psi(f) = f\circ p $ for $f\in\CC(X)$, we have that 
for all $f\in\CB$, 
$$
\Psi(f) = f\circ p \ \ (\rho\text{-almost 
everywhere}).
$$

We claim that the map $p\colon(W, \rho, R)\to (X, \mu, T)$ is an 
isomorphism between measure preserving systems.  
Indeed, 
the range of the map  $f\mapsto f\circ p \colon 
L^2(\mu)\to L^2(\rho)$ is closed in $L^2(\rho)$ 
because this map is 
an isometry, and it contains $\Psi(\CB) = \CC(W)$ 
and thus it is equal to $L^2(\rho)$.  In particular, $(W, \rho, R)$ 
is ergodic.  

Finally, for every function $g\in\CC(Y)$, we have that 
$g\circ q = \Psi(g\circ\pi) = g\circ\pi\circ p $ ($\rho$-almost everywhere)
and so 
$q = \pi\circ p$ ($\rho$-almost everywhere).

In particular, the image of $\rho$ under $q$ is $\nu$.  

\subsubsection{Construction of the sequence of intervals}  
Since $\rho$ is ergodic under $R$, it admits a generic point $w_1$. 
Recall that this means  that 
for every $f\in\CC(W)$, 
$$
\lim_{j\to+\infty}\frac{1}{j}\sum_{n=0}^{j-1}f(R^nw_1) = \int f\,d\rho \ .
$$
Set $x_1 = p(w_1)$.  Since $x_0$ is a transitive point of $X$, we can 
choose as in Section~\ref{subsubsec:averaging} 
a sequence of integers $(k_j\colon j\geq 1)$ such that
\begin{equation}
\label{eq:dX}
\lim_{j\to+\infty}\sup_{0\leq n\leq j}d_X(T^nx_1,T^{k_j+n}x_0)=0\ .
\end{equation}

Set $y_1 = q(w_1)$.  Let $\eta$ be a point 
in the closure of the sequence $(S^{k_j}\colon j\geq 1)$ 
in the Ellis semigroup~\cite{E} of $(Y,S)$.  Since $(Y,S)$ is distal, we have 
(see~\cite{auslander}, chapter $5$)
that $\eta$ is a bijection from $Y$ onto itself.  
Pick $y_0\in Y$ such that $\eta(y_0) = y_1$.  
Thus passing, if necessary, to a subsequence of $(k_j\colon j\geq 1)$, 
which we also denote by $(k_j\colon j\geq 1)$, we have that 
$T^{k_j}y_0$ converges to $y_1$.  Again replacing this sequence 
by a subsequence, we can assume that 
\begin{equation}
\label{eq:dY}
\lim_{j\to+\infty}\sup_{0\leq n< j}d_Y(S^ny_1, S^{k_j+n}y_0)= 0\ .
\end{equation}

For all $j\geq 1$, set $I_j = [k_j, k_j+j-1]$. 
Let $f$ be a continuous function on $X$ and $g$ a continuous function 
on $Y$. By~\eqref{eq:dX} and~\eqref{eq:dY} we have that 
\begin{gather*}
\lim_{j\to+\infty}\sup_{0\leq n < 
j}\bigl|f(T^nx_1)-f(T^{k_j+n}x_0)\bigr|=0
\text{ and }\\
\lim_{j\to+\infty}\sup_{0\leq n< j} 
\bigl|g(S^ny_1)-g(S^{k_j+n}y_0)\bigr|=0\ .
\end{gather*}
Thus 
\begin{equation}
\label{eq:f-g}
\lim_{j\to+\infty}\Bigl(\frac{1}{|I_j|}\sum_{n\in I_j}f(T^nx_0)g(S^ny_0) - 
\frac{1}{j}\sum_{n=0}^{j-1}f(T^nx_1)g(S^ny_1)\Bigr)= 0 \ .
\end{equation}
For each integer $n$, 
$$
f(T^nx_1)g(S^ny_1) = f\circ p(R^nw_1).g\circ q(R^nw_1)\ .
$$
Since $w_1$ is a generic point with respect to the measure $\rho$, 
the second average in~\eqref{eq:f-g} converges to 
$$
\int (f\circ p).(g\circ q)\,d\rho =
\int (f\circ p).(g\circ \pi\circ p) \,d\rho =
\int f.(g\circ \pi)\,d\mu 
$$
because $q= \pi\circ p$ ($\rho$-almost everywhere) and the image of $\rho$ 
under $p$ is $\mu$.
\end{proof}

\subsection{Inverse results} 

\begin{proposition}
\label{prop:extended2}
Let $k\geq 2$ be an integer, $\ua$ be a bounded sequence, and $\delta
>0$. Then there exists a sequence of intervals $\uI=(I_j\colon j\geq 
1)$ whose lengths tend to infinity and a $(k-1)$-step smooth nilsequence 
$\ub$ such that
\begin{enumerate}
\item
The sequence $\ua$ satisfies property $\CP(k)$ on $\uI$ and 
$\norm\ua_{\uI,k}\geq\norm\ua_{U(k)}-\delta$.
\item The sequence $\ua-\ub$ satisfies property $\CP(k)$ on $\uI$  
and
$\norm{\ua-\ub}_{\uI,k}<\delta$.
\end{enumerate}
\end{proposition}
\begin{proof}

 Let $(X,T,x_0)$ be the pointed dynamical system associated to 
the algebra spanned by the sequence $\ua$, 
as in Section~\ref{subsec:pointeddynamic}.  
Let $f$ be the continuous function on $X$ 
defined by $f(T^nx_0)=a_n$ for every $n\in\Z$.

By Proposition~\ref{prop:supergodic}, 
there exists an invariant ergodic 
measure $\mu$ on $X$ with $\nnorm f_{\mu,k}\geq 
\norm\ua_{U(k)}-\delta$.
By Corollary~\ref{cor:structure} of the Structure Theorem 
there exist a $(k-1)$-step 
nilsystem $(Y,S,\nu)$, a measure theoretic factor map $\pi\colon 
(X,\mu,T)\to (Y,\nu,S)$, and a smooth function $g$ on $Y$
 with $\nnorm{f-g\circ\pi}_{\mu,k}<\delta$.

Recall that every nilsystem is distal. 
Now, let $\uI$ and $y_0$ be given by Proposition~\ref{prop:extended} 
and let $\ub$ be the nilsequence given by $b_n=g(S^ny_0)$ for every 
$n\in\Z$.

 The measure on $X$ associated to $\uI$ as in~\ref{subsubsec:averaging} 
is equal to $\mu$.
Thus the sequence $\ua$ satisfies property $\CP(k)$ on $\uI$
and $\norm{\ua}_{\uI,k}=\nnorm 
f_{\mu,k}\geq\norm\ua_{U(k)}-\delta$. 
To prove 
Proposition~\ref{prop:extended2},
we are left with proving that the 
sequence $\ua-\ub$ 
satisfies property $\CP(k)$ on $\uI$ and that
$\norm{\ua-\ub}_{\uI,k}<\delta$.

For $h=(h_1,\dots,h_k)\in\Z^k$, we have 
\begin{multline*}
 \prod_{\epsilon\in \{0,1\}^k}C^{|\epsilon|}
(a_{n+\epsilon\cdot h}-b_{n+\epsilon\cdot h})
=  \prod_{\epsilon\in \{0,1\}^k}C^{|\epsilon|}
\bigl(f(T^{n+\epsilon\cdot h}x_0)-g(S^{n+\epsilon\cdot h}y_0)\bigr)\\
=\sum_{(A,B)\text{ partition of }\{0,1\}^k}
(-1)^{|B|}\prod_{\epsilon\in A}C^{|\epsilon|}
f(T^{n+\epsilon\cdot h}x_0)
\prod_{\epsilon\in B}C^{|\epsilon|}
g(S^{n+\epsilon\cdot h} y_0)\ .
\end{multline*}
By definition of $\uI$, the 
averages (with respect to $n$) of the above expression 
on this sequence of 
intervals converge to
\begin{multline*}
\sum_{(A,B)\text{ partition of }\{0,1\}^k}(-1)^{|B|}
\int \prod_{\epsilon\in A}C^{|\epsilon|}f(T^{\epsilon\cdot h}x)
\prod_{\epsilon\in B}C^{|\epsilon|}
g\circ\pi(T^{\epsilon\cdot h} x)\, d\mu(x)\\
=\int\prod_{\epsilon\in \{0,1\}^k}C^{|\epsilon|}
(f-g\circ\pi)(T^{\epsilon\cdot h}x)\,d\mu(x)\ .
\end{multline*}
By definition, the averages (with respect to $h\in\Z^d$) of the first term converge 
to $\norm{\ua-\ub}_{\uI,k}$ and, by Corollary~\ref{cor:dual}, the averages of the last integral converge to 
$\nnorm{f-g\circ\pi}_{\mu,k}<\delta$ and we are done.
\end{proof}

We now prove the Inverse Theorem (Theorem~\ref{th:inverse}).  We
recall the statement here for convenience. 
\begin{theorem*}
 Let $\ua = (a_n\colon n\in\Z)$ be a bounded sequence.  
 Then for every $\delta>0$, there exists a $(k-1)$-step smooth nilsequence 
$\ub=(b_n\colon n\in\Z)$  such that 
$$
 \nnorm \ub_k^*=1\text{ and }
\lim_{N\to+\infty}\,\sup_{M\in\Z}\Bigl|\frac{1}{N}\sum_{n=M}^{M+N-1}
a_nb_n\Bigr|
\geq \norm \ua_{U(k)}-\delta\ .
$$
\end{theorem*}

\begin{proof} We can assume without loss that $\norm\ua_{U(k)}>\delta$.
Let $\uI$ and $\uc$ be as in Proposition~\ref{prop:extended2}, but 
with $\delta/3$ instead of $\delta$; we 
write $c_n=g(S^ny_0)$ for $n\in\Z$, 
where $(Y,S,\nu)$ is an ergodic $(k-1)$-step 
nilsystem, $y_0\in Y$, and $g$ is a smooth function on $Y$. 
We define $h=\nnorm g_{k}^{-2^k+1}.\CD_kg$ and $\ub$ to be the 
sequence given by $b_n=h(S^ny_0)$,  and we check that the 
announced properties are satisfied. 

By Corollary~\ref{cor:smooth}, $h$ is a smooth function and 
$\nnorm h_{\nu,k}^*=1$ by~\eqref{eq:Dkfdualfunction} and thus 
$\nnorm\ub_k^*=1$. 
We have
\begin{multline*}
\LimAvI{\uI}{c_nb_n}=\LimAvI{\uI}{g(S^ny_0)h(S^ny_0)}=\int 
g.h\,d\nu=\nnorm g_k\\
=\norm\uc_{\uI,k}\geq \norm\ua_{\uI,k}-\delta/3\geq 
\norm\ua_{U(k)}-2\delta/3\ .
\end{multline*}
On the other hand, by the Direct Theorem~\ref{th:upperbound},
$$
 \LimsupAvI{\uI}{(a_n-c_n)b_n}\leq
\norm{\ua-\uc}_{\uI,k}\,\nnorm\ub_k^* \leq\delta/3
$$
and we conclude that the $\liminf$ of the averages on $\uI$ of 
$a_nb_n$ is $\geq \norm\ua_{U(k)}-\delta$ and we are done. 
\end{proof}

\subsection{Proof of Theorem~\ref{th:cond-conv}}
We recall the statement for convenience.
\begin{theorem*}
For a bounded sequence $\ua=(a_n\colon n\in\Z)$, 
the following are equivalent.
\begin{enumerate}
\item\label{it:unif}
For every $\delta >0$, the sequence $\ua$ can be written as
$\ua'+\ua''$, where $\ua'$ is a $(k-1)$-step nilsequence,
 and
$\norm{\ua''}_{U(k)}<\delta$.
\item\label{it:converges}
For every $(k-1)$-step nilsequence $\uc=(c_n\colon n\in\Z)$, 
the averages of $a_nc_n$ converge.
\end{enumerate}
\end{theorem*}


We recall that property~\eqref{it:converges} means that the 
averages
$$
 \frac 1{|I_j|}\sum_{n\in I_j}a_nc_n
$$
converge for every sequence of intervals $\uI=(I_j\colon n\geq 1)$ 
whose lengths tend to infinity. The common value of these limits is 
written $\LimAv{a_nc_n}$.
\begin{proof}
\eqref{it:unif} $\Longrightarrow$  \eqref{it:converges}
This implication is given by Corollary~\ref{cor:decomp-converges2}.

\medskip\noindent \eqref{it:converges} $\Longrightarrow$ 
\eqref{it:unif}  

Assume that the sequence $\ua$ satisfies~\eqref{it:converges}.
Let $\ub$ and $\uI$ be as in Proposition~\ref{prop:extended2},
 but with $\delta/3$ instead of $\delta$. 
Define $\ua'=\ub$ and we are left with showing that
$\norm{\ua-\ub}_{U(k)}<\delta$.

Assume that this does not hold. By Theorem~\ref{th:inverse}, 
there exists a $(k-1)$-step 
smooth nilsequence $\uc$  and a sequence of intervals 
$\uJ$ whose lengths tend to 
infinity with 
$$
\nnorm \uc_k^*=1\text{ and }\bigl|\LimAvI{\uJ}{(a_n -b_n)
c_n}\bigr|\geq2\delta/3\ .
$$
 Now, the sequence $(b_nc_n)$
 is a product of two $(k-1)$-step nilsequences and thus it is also 
a $(k-1)$-step nilsequence and its averages converge. By hypothesis, 
the averages of the sequence  $(a_nc_n)$ converge, and thus the 
averages of the sequence $(a_n -b_n)
c_n$ converge. Since $\uI$ and $\uJ$ are sequences of intervals 
whose lengths tend to infinity, 
\begin{align*}
 \bigl|\LimAvI{\uI}{(a_n-b_n)c_n}\bigr|= & \bigl|\LimAv{(a_n-b_n)c_n}|\\
= & \bigl|\LimAvI{\uJ}{(a_n-b_n)c_n}\bigr|\geq 2\delta/3\ .
\end{align*}

On the other hand, by the Direct Theorem 
(Theorem~\ref{th:upperbound})
$$
  \bigl|\LimAvI{\uI}{(a_n-b_n)c_n}\bigr|\leq 
\norm{\ua-\ub}_{\uI,k}\,\nnorm \uc_k^*<2\delta/3
$$
and we have a contradiction.
\end{proof}

\section{An application in ergodic theory}
\label{sec:application}

\subsection{Proof of Theorem~\ref{th:NilWW}}
\label{subsec:proofNilWW}
We now turn to the generalization of the Wiener-Wintner 
Ergodic Theorem, replacing the exponential sequence $e(nt)$
by an arbitrary nilsequence.
Throughout this Section, for each integer 
$N\geq 1$, we write $I_N$ for the interval 
$[0,N-1]$ and we let $\uI$ denote the sequence of intervals 
$(I_N\colon N\geq 1)$.

Let $(X,\mu,T)$ be an ergodic system, $\phi$ be a bounded measurable function on $X$, and 
fix an integer $k\geq 2$.  We build a subset $X_0$ 
of full measure of $X$ on which the conclusion of the Theorem holds for every $(k-1)$-step nilsequence $\ub$.

For every integer $r\geq 1$,  Corollary~\ref{cor:structure} of the Structure 
Theorem provides a $(k-1)$-step nilsystem $(Z_r,\nu_r,S_r)$, a factor 
map $\pi_r\colon X\to Z_r$ and a continuous function $f_r$ on $Z_r$ 
such that 
$$\nnorm{\phi-f_r\circ\pi_r}_k<r\inv\ .$$

By Corollary~\ref{cor:normergo}, there exists a 
subset $E_r$ of $X$ with $\mu(E_r)=1$ such that for every $x\in 
E_r$, we have
$$
 \norm{(\phi(T^nx)-f_r\circ\pi_r(T^nx)\colon 
n\in\Z)}_{\uI,k}=\nnorm{\phi-f_r\circ\pi_r}_k\leq r\inv\ .
$$
Note that we consider the map $\pi_r$ to be defined everywhere. For 
$\mu$-almost every $x$, we have that 
$f_r\circ\pi_r(T^nx)=f_r(S_r^n\pi_r(x))$ for every $n\in\Z$.
Therefore, there exists a set $E'_r\subset X$ with $\mu(E'_r)=1$ 
such that 
$$
\norm{(\phi(T^nx)-f_r(S_r^n\pi_r(x))\colon 
n\in\Z)}_{\uI,k}=\nnorm{\phi-f_r\circ\pi_r}_k\leq r\inv
$$
for every $x\in E'_r$. 

Set $X_0 = \bigcap_{r=1}^\infty E'_r$.  
For every $x\in X_0$, the sequence $(\phi(T^nx)\colon n\in\Z)$ 
satisfies the hypothesis of Corollary~\ref{cor:decomp-converges}, 
completing the proof. \qed

\subsubsection{Proof of Corollary~\ref{cor:WW}}
 Let  $(X,\mu,T)$ be an ergodic system, $\phi$ be a bounded measurable 
function on $X$, and let $X_0$ be the subset of $X$ introduced in 
Theorem~\ref{th:NilWW}. 
Let $x\in X_0$ and $p$ be a generalized polynomial.

For every $n\in\Z$, let $\{p(n)\}$ denote the fractional part of $p(n)$. Then 
$\{p(\cdot)\}$ is a bounded generalized polynomial. 
In~\cite{BL} (Theorem A, (ii)), 
it is shown that there exist  an ergodic 
 nilsystem $(Y,\nu,S)$, a point $y\in Y$, and a 
Riemann integrable function $f$ on $Y$ with $\{p(n)\}=f(S^ny)$ for 
every $n\in\Z$.

For 
every $\delta>0$, there exists a continuous function $g$ on $Y$ with 
$\norm{f-g}_{L^1(\nu)}\leq\delta$. The sequence $(g(S^ny)\colon 
n\in\Z)$ is a nilsequence and thus by definition of $X_0$, the averages
on $\uI$ of $\phi(T^nx)g(S^ny)$ converge. On the  other hand, since 
the function $|f-g|$ is Riemann integrable and $(Y,S)$ is uniquely 
ergodic, we have that
\begin{multline*}
\LimsupAvI{\uI}{\phi(T^nx)(f(S^ny)-g(S^ny))}\\
\leq\norm\phi_\infty\,
 \LimAvI{\uI}{|f(S^ny)-g(S^ny)|}=
\norm\phi_\infty\int|f-g|\,d\nu\leq\norm\phi_\infty\delta\ .
\end{multline*}
Therefore the averages on $\uI$ of $\phi(T^nx)\{p(n)\}=\phi(T^nx)f(S^ny)$ form a Cauchy 
sequence.

 We remark that for every $n\in\Z$, we have that
$e(p(n))=e(\{p(n)\})=e(f(S^ny))$ and that the function $e(f(\cdot))$ 
is Riemann integrable on $Y$.
The same proof gives the second claim of the corollary.
\qed

\subsection{Examples}
\label{subsec:morse}

Similar methods can be used to show 
show that some explicit sequences satisfy 
the hypothesis~\eqref{it:unif} of Theorem~\ref{th:cond-conv} 
and thus are universally 
good for the convergence in norm of multiple ergodic averages.

\begin{proposition}
\label{prop:riemann2}
Let $(X,T)$ be a uniquely ergodic system with invariant measure $\mu$ 
and let $k\geq 2$ be an integer. Let $(Z_{k},\mu_{k},T)$ be the factor 
defined in the Structure Theorem (Theorem~\ref{th:structure}) and 
assume that the factor map $\pi_{k}\colon X\to Z_{k-1}$ is continuous.
Let $f$ be a Riemann integrable function on $X$ and let $x\in X$. Then 
the sequence $(f(T^nx)\colon n\in\Z)$ satisfies
hypothesis~\eqref{it:unif} of Theorem~\ref{th:cond-conv}.
\end{proposition}

\begin{proof} Let $\ua$ be the sequence  $(f(T^nx)\colon n\in\Z)$ and 
let $\delta>0$. We want to show that we can write $\ua=\ua'+\ua''$ where 
$\ua'$ us a $(k-1)$-step nilsequence and $\norm{\ua''}_{U(k)}<\delta$.

Let $(Y,S,\nu)$, $p\colon X\to Y$, and $h$ be the $(k-1)$-step 
nilsystem, the factor map, and the function on $Y$ given by 
Corollary~\ref{cor:structure}. Recall that $Z_k$ is the inverse limit 
(in both the topological and measure theoretical senses) of all factors of 
$X$ which are $(k-1)$-step nilsystems~\cite{HK}. Thus $Y$ is a factor of $Z_k$ 
and the factor map $q\colon Z_k\to Y$ is continuous.
Therefore the factor map $p=q\circ\pi_k$ mapping $X\to Y$ is continuous. 

We define the sequences $\ua'$ and $\ua''$ by $a'_n=h\circ p(T^n x)$ 
 and $a''_n=f(T^nx)-h\circ p(T^nx)$ for every 
$n\in\Z$.
 Then $\ua'$ is a $(k-1)$-step nilsequence. Since the function $h\circ 
p$ is continuous, the function $f-h\circ 
p$ is Riemann integrable, and Corollary~\ref{cor:normtotergo} implies that
$\norm{\ua''}_{U(k)}=\nnorm{f-h\circ p}_k<\delta$.
\end{proof}

We use this proposition to prove Proposition~\ref{prop:genpol} on
generalized polynomials.
\begin{proof}[Proof of Proposition~\ref{prop:genpol}]
Let $p$ be a generalized polynomial.  For every $n\in\Z$, let $\{p(n)\}$ 
denote the fractional part of $p(n)$. We begin with the same argument as in 
the proof of Corollary~\ref{cor:WW}.

There exists an integer $\ell\geq 1$, an ergodic 
$\ell$-step nilsystem $(X=G/\Gamma,\mu,T)$, a point $x\in X$, and a 
Riemann integrable function $f$ on $X$ with $\{p(n)\}=f(T^nx)$ and 
$e(p(n))=e(\{p(n)\})=e(f(T^nx))$ 
for every $n\in\Z$.

The system $(X,\mu,T)$ satisfies the hypotheses of 
Proposition~\ref{prop:riemann2}. Indeed, for $k>\ell$ we have that $Z_k=X$ 
and for $k<\ell$, $Z_k$ is the quotient $G/G_k\Gamma$ of $X$. The 
result follows.
\end{proof}

We now prove Proposition~\ref{prop:morse}, which states that the 
Thue-Morse sequence satisfies also the hypothesis of Theorem~\ref{th:cond-conv}.

\begin{proof}[Proof of Proposition~\ref{prop:morse}]

Let $\ua = (a_n\colon n\in\Z)$ be the Thue-Morse sequence. We recall 
some of its properties (see~\cite{Q}). 

There exists a uniquely ergodic system $(X,T,\mu)$, a point $x_0\in 
X$, and a continuous function $\phi$ on $X$ with $a_n=\phi(T^nx_0)$ for 
every $n\in\Z$. Moreover, the factor map $\pi_1\colon X\to Z_1$ on the 
Kronecker factor $Z_1$ of $X$ is continuous. Finally, the map $\pi$ 
is two to one almost everywhere.

For every integer $k\geq 2$, the factor $Z_k$ of $X$, as given by
the Structure Theorem, is an extension of 
$Z_{k-1}$ by a connected compact abelian group~\cite{HK}. It follows 
that $Z_k=Z_1$ for every $k$.

Therefore the hypotheses of Proposition~\ref{prop:riemann2} are 
satisfied and we are done.
\end{proof}

\subsection{Proof of Theorem~\ref{th:convnilconvmult}}
We now prove the generalization of the spectral theorem.  
Starting with an arbitrary measure preserving 
system $(Y,S,\nu)$, by ergodic decomposition we can assume 
that $(Y,S,\nu)$ is an ergodic system.

We recall the following result from~\cite{HK} (Theorem 12.1):
\begin{theorem*}
Let $g_0,\dots, g_{k-1}$ be measurable functions on $(Y,S,\nu)$ with 
$\norm{g_i}_\infty\leq 1$ for $i\in\{0, \ldots, k-1\}$.
Then
$$
 \limsup_{N\to+\infty}
 \Bigl|\frac 1{N}\sum_{n=0}^{N-1}\int \prod _{i=0}^{k-1}
S^{in}g_i\,d\nu\Bigr|\leq c\min_{i\in\{0,\ldots, k-1\}}\nnorm{g_i}_{k-1}
$$
where $c$ is a constant depending only on $k$.
\end{theorem*}

Proceeding as in~\cite{BHK} (proof of Corollary 4.5 from Theorem 4.4), we 
deduce:
\begin{corollary}
\label{cor:linear1}
Let $g_0,\dots, g_{k-1}$ be measurable functions on $(Y,S,\nu)$ with 
$\norm{g_i}_\infty\leq 1$ for $i\in\{0,\ldots, k-1\}$. 
Then 
$$
 \limsup_{N\to+\infty}\frac 1{N}\sum_{n=0}^{N-1}
\Bigl| \int \prod _{i=0}^{k-1} 
S^{in}g_i\,d\nu\Bigr|^2\leq c^2\min_{i\in\{0,\ldots, k-1\}}\nnorm{g_i}_{k}^2\ .
$$
\end{corollary}

We deduce:

\begin{corollary}
\label{cor:linear2}
Let $f_1,\dots,f_k$ be bounded functions on $(Y,S,\nu)$ with 
$\norm{f_i}_\infty\leq 1$ for $i\in\{1,\ldots,k\}$ and let $\ua=(a_n\colon n\in\Z)$ 
be a sequence with $\norm\ua_\infty\leq 1$. Then
\begin{equation}
\label{eq:linear2}
 \limsup_{N\to+\infty}\Bigl\Vert\frac 1{N}\sum_{n=0}^{N-1}
a_n\prod_{i=1}^kS^{in}f_i\Bigr\Vert_{L^2(\nu)}\leq 
k^{1/4} c^{1/2}\min_{i\in\{1,\ldots, k\}}\nnorm{f_i}_{k+1}\ .
\end{equation}
\end{corollary}

\begin{proof} Let $\ell\in\{1,\dots,k\}$ be such that 
$\norm{f_\ell}_{k+1}=\min_{i\in\{1, \ldots, k\}}\nnorm{f_i}_{k+1}$ and
let $Q$ be the $\limsup$ in the left hand side of~\eqref{eq:linear2}.

By the van der Corput Lemma (Appendix~\ref{ap:vdc}):
$$
 Q^2\leq \limsup_{M\to+\infty}\frac 1M\sum_{m=0}^{M-1}\Bigl|
\limsup_{N\to+\infty}\frac 1{N}\sum_{n=0}^{N-1}
\overline{a_n}a_{n+m}
\int\prod _{i=1}^k S^{in}(\overline{f_i}.S^{im}f_i)\,d\nu\Bigr|\ .
$$
By the Cauchy-Schwarz Inequality,
\begin{multline*}
 Q^4\leq \limsup_{M\to+\infty}\frac 1M\sum_{m=0}^{M-1}
\limsup_{N\to+\infty}\frac 1{N}\sum_{n=0}^{N-1}
\Bigl| \int \prod _{i=1}^k 
S^{in}(\overline{f_i}.S^{im}f_i)\,d\nu\Bigr|^2\\
=\limsup_{M\to+\infty}\frac 1M\sum_{m=0}^{M-1}
\limsup_{N\to+\infty}\frac 1{N}\sum_{n=0}^{N-1}
\Bigl| \int \prod_{i=0}^{k-1} 
S^{in}(\overline{f_{i+1}}.S^{(i+1)m}f_{i+1})\,d\nu\Bigr|^2\ .
\end{multline*}
Applying Corollary~\ref{cor:linear1} 
to the functions $g_i=\overline{f_{i+1}}.S^{(i+1)m}f_{i+1}$, we have that
\begin{multline*}
 Q^4\leq c^2\limsup_{M\to+\infty}\frac 1M\sum_{m=0}^{M-1}
\nnorm{\overline{f_\ell}.S^{\ell m}f_\ell}_k^2
\leq k c^2
\limsup_{M\to+\infty}\frac 1{kM}\sum_{m=0}^{kM-1}
\nnorm{\overline{f_\ell}.S^{ m}f_\ell}_k^2\\
\leq k c^2 \Bigl(\limsup_{M\to+\infty}\frac 1{kM}\sum_{m=0}^{kM-1}
\nnorm{\overline{f_\ell}.S^{ m}f_\ell}_k^{2^k}\Bigr)^{1/2^{k-1}}
\end{multline*}
by the H\"older Inequality. 
 By~\eqref{eq:k-to-kplus1}, the last $\limsup$ is actually 
a limit and is equal to $\nnorm{f_\ell}_{k+1}^4$ and we are done.
\end{proof}

We now return to the proof of Theorem~\ref{th:convnilconvmult}. 
We assume that $\ua=(a_n\colon n\in\Z)$ is a bounded sequence 
such that the averages 
$$
\frac{1}{N}\sum_{n=0}^{N-1}a_nb_n
$$
converge as $N\to+\infty$ 
for every $k$-step nilsequence $\ub=(b_n\colon n\in\Z)$.  
We assume that $(Y,S,\nu)$ is an ergodic system and 
$f_1, \ldots, f_k\in L^\infty(\nu)$.  We show the convergence 
of the averages
$$
\frac{1}{N}\sum_{n=0}^{N-1}a_nS^nf_1\ldots S^{kn}f_k
$$
in $L^2(\nu)$.  

Let $Z_k$ be the $k$-th factor of $(Y,S,\nu)$, 
as given by the Structure Theorem.  If for some $i\in\{1, \ldots, k\}$ we 
have $\E(f_i\mid Z_k)=0$, then  $\nnorm{f_i}_{k+1}=0$.  Then by 
Corollary~\ref{cor:linear2}, the above 
averages converge to zero in 
$L^2(\nu)$. We say that the {\em factor $Z_k$ is characteristic for 
the convergence of these averages}.

Therefore, in order to prove the convergence of these 
averages, for arbitrary bounded functions, it 
suffices to prove the convergence when the functions are measurable 
with respect to the factor $Z_k$.

By the Structure Theorem, $Z_k$ is an inverse limit of $k$ step 
nilsystem.  Thus by density, we can assume that the
functions $f_i$ are measurable with respect to a $k$-step nilsystem 
$(Z,S)$ which 
is a factor of $(Y,S,\nu)$. 
By density again, we are reduced to the case that $(Y,\nu,S)$ is a $k$-step nilsystem 
and that the functions  $f_1,\dots,f_k$ are continuous.

But in this case, for every $y\in Y $ the sequence 
$$
(f_1(S^ny).f_2(S^{2n}y).\cdots.f_k(S^{kn}y)\colon n\in\Z)
$$
 is a $k$-step nilsequence and by hypothesis, the averages
$$
\frac 
1N\sum_{n=0}^{N-1}a_n\,f_1(S^ny).f_2(S^{2n}y).\cdots.f_k(S^{kn}y)
$$
converge for every $y\in Y$.\qed

\appendix

\section{The van der Corput Lemma}
\label{ap:vdc}
We state the van der Corput Lemma, as used in our 
set up (see~\cite{KN}):
\begin{vdC}
Let $\ua = (a_n\colon n\in\Z)$ be a sequence with $|a_n|\leq 1$ 
for all $n\in\Z$ and let $I$ be an interval in $\Z$.  Then for 
every integer $H\geq 1$, we have 
$$
\vert \frac{1}{|I|}\sum_{n\in I}a_n\vert^2 \leq \frac{4H}{|I|}+
\Bigl|\sum_{h=-H}^H\frac{H-|h|}{H^2}\;\frac{1}{|I|}
\sum_{n\in I}a_{n+h}\overline{a_n}\Bigr|
\ .
$$
\end{vdC}

\section{Parallelepipeds in nilmanifolds}
\label{appendix:cubes}

We explain the cubic structure associated to a nilmanifold.
In the literature, 
there are (at least) two presentations of these objections, 
in~\cite{HK} and in Appendix E of~\cite{GT}.
The results proved in these papers are often  recalled here 
without proof, but we need a bit more than just those results. 
We use the notation of~\cite{HK}. 
The group that we denote by $G\type k_{k-1}$ is the same as the 
group $\mathop{\HP}^k$ of~\cite{GT}.  

The $k$'s in index and exponent that occur everywhere are cumbersome 
but necessary as we use an induction at some point.

\subsection{Algebraic preliminaries}

We begin with some algebraic constructions 
involving ``cubes.''
Let $G$ be a group and $k\geq 1$ be an integer.

\subsubsection{Two constructions of the ``side group''}

We use the notation of Section~\ref{sec:ergseminorms}.
We write
$\bzero=(0,0,\dots,0)\in\{0,1\}^k$ and $\one=(1,1,\dots,1)\in\{0,1\}^k$.

As before, if $X$ is a set, $X\type k=X^{2^k}$ and points of $X\type k$ are 
written as $\ulx=(x_\epsilon\colon 
\epsilon\in\{0,1\}^k)$.
For $x\in X$, $x\type k\in X\type k$ is the element
$(x,x,\dots,x)$, with $x$ repeated $2^k$ times.
If $f\colon X\to Y$ is a map, $f\type k\colon X\type k\to Y\type k$ 
denotes the diagonal map: $(f(\ulx))_\epsilon=f(x_\epsilon)$ for all
$\epsilon\in\{0,1\}^k$.

For $g\in G$ and $1\leq i\leq k$, $g_i\type k=((g_i\type k)_\epsilon\colon
\epsilon\in\{0,1\}^k)$ is given by:
$$
 \bigl(g\type k_i)_\epsilon=
\begin{cases} g & \text{if }\epsilon_i=1 \\
 1 &\text{if }\epsilon_ i=0\ .
\end{cases}
$$
(Note that we mean $\epsilon = (\epsilon_1, \ldots, 
\epsilon_{k})$.)
$G\type k_{k-1}$ is the subgroup of $G\type k$  spanned by 
$$
 \{ g\type k\colon g\in G\}\cup\{ g_i\type k\colon 1\leq i\leq k,\ 
g\in G\}\ .
$$

The same group was also introduced in~\cite{GT}, but with 
a different definition and notation. 
We recall their presentation, but in our notation, substituting 
``upper faces'' for ``lower faces'' for coherence.
We start with some notation.

It is convenient to view $\{0,1\}^k$ as the 
set of vertices of the unit Euclidean cube.
If $J$ is a subset of $\{ 1,\dots,k \}$ and $\eta\in \{0,1\}^J $, 
the set 
$$
\alpha=\{ \epsilon\in\{0,1\}^k\colon \epsilon_i=\eta_i\text{ for all }i\in J\}
$$
is called a \emph{face} of $\{0,1\}^k$.  
The dimension of $\alpha$ is $\dim(\alpha)=k-|J|$.
If all coordinates of 
$\eta$ are equal to $1$, then this face is called an \emph{upper face}. 
In particular, $\alpha_0=\{0,1\}^k$ is the unique upper face of dimension $k$, 
corresponding to $J=\emptyset$; $\{\one\}$ is the unique upper face of 
dimension zero, corresponding to $J=\{1,\dots,k\}$. 
The $k$ upper faces of dimension $k-1$ are 
$\alpha_i=\{\epsilon\in\{0,1\}^k\colon\epsilon_i=1\}$ for
 $1\leq i\leq k$. 
Let  $\alpha_0,\alpha_1,\dots,\alpha_{2^k}$ be an enumeration of all of 
the upper faces such that $\alpha_0,\dots,\alpha_k$ are as above and 
$\dim(\alpha_i)$ is a decreasing sequence; in particular, 
$\alpha_{2^k}=\{\one\}$. 

If $\alpha$ is a face and $g\in G$, we write $g_\alpha\type
k=((g_\alpha\type k)_\epsilon\colon \epsilon\in\{0,1\}^k)$ for the element
of $G\type k$ given by:
$$
 \bigl(g\type k_\alpha\bigr)_\epsilon=\begin{cases}
g&\text{if }\epsilon\in\alpha\ ;\\ 
1& \text{otherwise\ .}
\end{cases}
$$
In particular, the elements $g\type k_i$ defined above can be written 
as $g\type k_{\alpha_i}$.

In~\cite{GT}, $\HP^k(G)$ is defined 
to be the set of elements $\ulg\in G\type k$ that 
can be written as 
\begin{equation}
\label{eq:def2group}
 \ulg=(g_1)\type k_{\alpha_1}(g_2)\type k_{\alpha_2}\dots 
(g_{2^k})\type k_{\alpha_{2^k}}
\text{ where }g_i\in G_{k-\dim(\alpha_i)}\text{ for every }i\in\{1, \ldots, k\}\ .
\end{equation}
Here $G_0=G_1=G$; 
in all other places in the paper, we use $G_0$ to denote a
different object (the connected component of the identity of $G$).

Let us explain briefly why $G\type k_{k-1}$ and $\HP^k(G)$ are 
actually equal.
By a direct computation, Green and Tao show that $\HP^k(G)$ is a subgroup of 
$G\type k$; since it contains the generators of $G\type k_{k-1}$, it 
contains this group. On the other hand, it is shown in~\cite{HK} (section 5)
that for every side $\alpha$ of dimension $d$ and every 
$g\in G_{k-\dim(\alpha)}$, $g\type k_\alpha$ belongs to $G\type 
k_{k-1}$ (and more precisely to $(G\type k_{k-1})_{k-\dim(\alpha)}$) 
and thus $\HP^k(G)\subset G\type k_{k-1}$.  We have equality. 

In the sequel we only use the notation $G\type k_{k-1}$. Depending on
the property to be proven, the first or second presentation 
is more convenient.

\subsubsection{Algebraic properties}
We have:
\begin{enumerate}
\item\label{it:coordGamma}
Let $\Gamma$ be a subgroup of $G$. 
If all coordinates of $\ulg$  belong to 
$\Gamma$ except possibly $g_\bzero$, then $g_\bzero\in\Gamma G_k$.
\item\label{it:coordG}
In particular, if all coordinates of $\ulg\in G\type k_{k-1}$ are 
equal to $1$ except possibly $g_\bzero$, then $g_\bzero\in G_k$.
\end{enumerate}
\trucenumi

The second statement is proved (in a perhaps concealed place) in~\cite{HK}
via induction on 
$k$, and the first one is not stated explicitly but follows with  
a similar proof. Both statements follow easily from the second 
definition  of $G\type k_{k-1}$ and the symmetry of this set, 
allowing us to substitute the coordinate $g_\one$ for $g_\bzero$.

We need two more groups for our proofs. In this appendix, we write
$$
H_k=\{\ulg\in G\type k_{k-1}\colon g_\bzero=1\}
\text{ and }
G\type k_k=\{g\type k\colon g\in G\}\ .
$$
(The first group is not defined in the papers.)
Then $H_k$ is clearly a normal subgroup of $G\type 
k_{k-1}$ and $G\type k_{k-1}=H_k.G\type k_k$.
Moreover, $H_k$ is the group spanned by the elements 
$g\type k_i$ for $1\leq i\leq k$ and $g\in G$; in the second 
presentation of $G\type k_{k-1}$, it consists of elements that can be 
written as in~\eqref{eq:def2group} with $g_1=1$.

We have
\begin{enumerate}
\enumitruc
\item
\label{it:K2}
$(H_k)_2=H_k\cap (G_2)\type k$.
\item
\label{it:G2}
$(G\type k_{k-1})_2=G\type k_{k-1}\cap (G_2)\type k$.
\end{enumerate}
\trucenumi

\begin{proof} 
We prove~\eqref{it:K2}.
The inclusion $(H_k)_2\subset H_k\cap (G_2)\type k$ is obvious.

Let $\alpha$ be a face of dimension $d<k-1$ containing $\one$.
Let $g\in G$ and $h\in G_{k-d-1}$. We can chose a face $\beta$ of 
dimension $k-1$ and a face 
$\gamma$ of dimension $d+1$ such that $\alpha=\beta\cap\gamma$. We have
$$
g\type k_\alpha\in H_k\ ;\ h\type k_\gamma\in H_k
\text{ and }
 [g;h]\type k_\alpha=\bigl[ g\type k_\beta;h\type k_\gamma\bigr]\ .
$$
Thus $[g;h]\type k_\alpha\in (H_k)_2$. 
Therefore, for any $q\in G_{k-d}$, we have that $q\type k_\alpha\in 
(H_k)_2$. 

Using this remark, we can show the inclusion 
$H_k\cap (G_2)\type k\subset (H_k)_2$. 
Let $\ulg$ be in the first of these groups. We write $\ulg$  as 
in~\eqref{eq:def2group} with $g_1=1$. By the remark, all terms 
of the form $(g_j)\type k_{\alpha_j}$ with $\dim(\alpha_j)<k-1$ in 
the product belong to  $(H_k)_2$ and we are reduced to 
show that the product of the $k$ remaining terms also belongs to this group.
We remark that all coordinates of this product belong to $G_2$.

Let $g\type k_\alpha$ be one of these terms. Then $\alpha$ is an upper face 
of dimension $k-1$ and it is immediate that there 
 exists $ \eta\in\{0,1\}^k$ such that 
$\eta$ belongs to $\alpha$ and does not belong to any other upper face of 
dimension $k-1$.
Therefore, $g$ is the coordinate $\eta$ of the product and $g\in 
G_2$. It follows that $g\type k_\alpha$ belongs to $(H_k)_2$ and we are done.

We now deduce~\eqref{it:G2}. Again, the inclusion 
$(G\type k_{k-1})_2\subset G\type k_{k-1}\cap (G_2)\type k$ is obvious.
Let $\ulg\in G\type k_{k-1}\cap (G_2)\type k$. We write  $\ulg=h\type 
k\ulq$ where $h\in G$ and $\ulq\in H_k$. We have that $g_0=h$ and so $h\in 
G_2$.  Thus $h\type k\in (G_2)\type k$. Moreover, $\ulq\in 
H_k\cap(G_2)\type k$ and by the second part of the Lemma, $\ulq\in 
(H_k)_2\subset (G_2)\type k_{k-1}$.
\end{proof}

\subsection{Topological properties}

Henceforth $G$ is a $r$-step nilpotent Lie group, $\Gamma$ is a discrete 
cocompact subgroup, and $X=G/\Gamma$. In applications $r$ will be 
equal to $k-1$ but the general case is used in an induction below.

In~\cite{HK} and~\cite{GT}, it is shown that 
\begin{enumerate}
\enumitruc
\item\label{it:Gkclosed}
$G\type k_{k-1}$ is a closed subgroup of 
$G\type k$ and hence is an $r$-step nilpotent Lie group.
\item
The group 
$\Lambda_k:=\Gamma\type k\cap G\type k_{k-1}$ is a  cocompact subgroup 
of $G\type k_{k-1}$.
\end{enumerate}
\trucenumi
 We do not reproduce the proof here. We define:
$$X_k= G\type k_{k-1}/(\Gamma\type k\cap G\type k_{k-1})\ .$$
For the moment we write 
$\nu_k$ for the Haar measure of $X_k$. 

The image of 
$\nu_k$ under the  projection $\ulx\mapsto x_\bzero$ is equal to the 
Haar measure $\mu$ of $X$.
We have that:
\begin{enumerate}
\enumitruc
\item
The group $\Theta_k:=H_k\cap \Gamma\type k$ is cocompact in $H_k$.
\end{enumerate}
\trucenumi

\begin{proof}
Every $\ulg\in H_k$ belongs to $G\type k_{k-1}$ and thus is at a bounded 
distance from some $\ulgamma\in \Lambda_k$.   Since $g_\bzero=1$, 
$\gamma_\bzero$ is at 
a bounded distance from $1$.  Since $\Gamma$ is discrete,
$\gamma_\bzero$ belongs to a finite subset $F$ of $\Gamma$.
 
We have that $\ulg$ is at a bounded distance from
 $((\gamma_\bzero)\type k)\inv\ulgamma$, which belongs to 
$G\type k_{k-1}\cap H_k=\Theta_k$.
\end{proof}

We define $W_k=H_k/\Theta_k\ $.  
Then $W_k$ is a $(k-1)$-step nilmanifold, naturally included in $X_k$
as a closed subset.

For every $g\in G$ we have that $g\type k$ belongs to $G\type 
k_{k-1}$. We deduce that for every $x\in X$, we have that $x\type 
k:=(x,x,\dots,x)$ belongs to $X_k$.

For every $x\in X$, we write
$$
W_{k,x}=\{ \ulx\in X_k\colon x_\bzero=x\}\ .
$$
We show:
\begin{enumerate}
\enumitruc
\item\label{it:WxHk}
Let  $x\in X$ and $g$ be a lift of $x$ in $G$. 
Then $W_{k,x}=g\type k. W_k$.
\end{enumerate}
\trucenumi
\begin{proof}
Let $\ulx\in W_{k,x}$ and $\ulh$ be a lift of $\ulx$ in $G\type k_{k-1}$.
Since $x_\bzero=x$, we have that $h_\bzero=g\gamma$ for some 
$\gamma\in\Gamma$. Let $\ulq=(g\type k)\inv\ulh(\gamma\type k)\inv$.
Then $\ulq\in H_k$ and its image $\uly$ in $H_k$ satisfies 
$g\type k\uly=\ulx$. We thus have that $W_{k,x}\subset g\type k.W_k$ and 
the opposite inclusion is obvious.
\end{proof}

\subsection{Dynamical properties}
Henceforth, we assume that $X$ is endowed with the translation $T$ by 
some $\tau\in G$ and that $(X,T,\mu)$ is ergodic. Recall that the 
same nilmanifold can be represented as a quotient in different ways.
As usual we assume that $G$ is spanned by the connected component 
$G_0$ of the identity  and  $\tau$. We claim that:

\begin{enumerate}
\enumitruc
\item\label{it:G0}
 $(G\type k_{k-1})_0=(G_0)\type k_{k-1}$. 
\item
\label{it:GG0}
$G\type k_{k-1}$ is spanned by $(G\type k_{k-1})_0$, $\tau\type k$, 
and the elements $\tau\type k_i$, $1\leq i\leq k$.
\item
\label{it:Hkconnected}
$H_k$ is spanned by $(H_k)_0$ and the elements $\tau\type k_i$, $1\leq 
i\leq k$.
\end{enumerate}
\trucenumi
\begin{proof}

By hypothesis and the first definition of $G\type k_{k-1}$, this group is 
spanned by elements of the form $g\type k$ for $g\in G_0$, $g\type 
k_i$  for $g\in G_0$ and $1\leq i\leq k$, $\tau\type k_i$ for $1\leq 
i\leq k$ and $\tau\type k$. This proves~\eqref{it:GG0}.

The commutator of two elements of the above type belongs to 
$(G_2)\type k_{k-1}\subset (G_0)\type k_{k-1}$, because it follows from
our assumption that $G_2\subset G_0$.
Then every element $\ulg$ of $G\type k_{k-1}$ can be written  as
$\ulg=\ulh(\tau\type 
k)^n(\tau\type k_1)^{m_1}\dots(\tau\type k_k)^{m_k}$ 
with $\ulh\in (G_0)\type k_{k-1}$.

If $\ulg\in (G\type k_{k-1})_0$, then by looking at the coordinate 
$\bzero$ of $\ulg$ we have that $h_0\tau^n= g_0$ belongs to 
$G_0$.  Thus $\tau^n\in G_0$. 

Let $i\in\{1,\dots,k\}$. As in the proof of~\eqref{it:K2}, there exists 
$\eta\in\{0,1\}^k$ such that $\tau_i\type k=\tau$ and $\tau_j\type 
k=1$ for $j\neq i$. We have that $g_\eta=h_\eta\tau_i^{m_i}$ and thus 
$\tau^{m_i}\in G_0$.  Thus $(\tau\type k_i)^{m_i}\in (G_0)\type k_{k-1}$.
This achieves the proof of~\eqref{it:G0}.

Now assume that $\ulg\in (H_k)_0$. Then it belongs to $(G\type 
k_{k-1})_0$ and we write it as above,
$\ulg=\ulh(\tau\type k_1)^{m_1}\dots(\tau\type k_k)^{m_k}$ with $\ulh\in 
(G_0)\type k_{k-1}$. We have that $h_\bzero=g_\bzero=1$ and so
$\ulh\in H_k\cap (G_0)\type k_{k-1}$ and this element belongs to 
$(H_k)_0$. This proves~\eqref{it:Hkconnected}.
\end{proof}

\begin{enumerate}
\enumitruc
\item
$X_k$ is ergodic under the action of $T\type k$ and $T\type k_i$, 
$1\leq i\leq k$.
\item
\label{it:Wkerg}
$W_k$ is ergodic under the transformations $T\type k_i$, 
$1\leq i\leq k$.
\end{enumerate}
\trucenumi

\begin{proof}
Let $Z$ be the compact abelian group $G/\Gamma G_2$ and $\sigma$ be the 
image of $\tau$ in $Z$. Since $T$ is ergodic, the translation by 
$\sigma$ on $Z$ is ergodic.

By~\eqref{it:G2} and (any) definition of $G\type k_{k-1}$,
 the quotient $G\type k_{k-1}/(G\type 
k_{k-1})_2\Lambda_k$
can be identified with the subgroup $Z\type k_{k-1}$ of $Z\type k$.
This group consists of the points $\ulz$ of $Z\type k$ which can be 
written as
$$
 \ulz=\bigl( u\prod_{i=1}^k v_i^{\epsilon_i}\colon 
\epsilon\in\{0,1\}^k\bigr)
$$
for some $ u,v_1,\dots,v_k\in Z$.
The transformations induced on this group by the transformations 
$T\type k$ and $T\type k_i$, $1\leq i\leq k$, are the translations by 
$\sigma\type k$ and $\sigma\type k_i$. In the above parametrization 
of $Z\type k_{k-1}$, these transformations correspond to the map 
$u\mapsto \sigma u$ and to the maps $v_i\mapsto \sigma v_i$, 
respectively.

Since the translation by $\sigma$ on $Z$ is ergodic, it follows easily 
that $Z\type k_{k-1}$ is ergodic under the translations 
by $\sigma\type k$ and $\sigma\type k_i$.
By~\eqref{it:G0} and Theorem~\ref{th:multiergo}, $X_k$ is ergodic under the 
action of $T\type k$ and $T\type k_i$, $1\leq i\leq k$.

The second statement is proved in the same way.
\end{proof}

We show:
\begin{enumerate}
\enumitruc
\item\label{it:Haar}
The Haar measure $\nu_k$ of $X_k$ is equal to the measure $\mu\type k$ 
defined in~\cite{HK} and described in Section~\ref{sec:ergseminorms}.
\end{enumerate}
\trucenumi
This result is proved in~\cite{HK}, but the context is so different from the
present one that we prefer to give a complete proof here.
\begin{proof}
We use induction on $k$.
By definition, $G\type 2_1=G\times G$ and so $X_1=X\times X$ and 
$\nu_1=\mu\times\mu$, which is equal to the measure $\mu_1$ of~\cite{HK}.

Assume that the announced property holds up to $k-1$ for some $k>1$.
In order to show the property for $k$, it suffices to show that 
when $f_\epsilon$, $\epsilon\in\{0,1\}^k$, are $2^k$ continuous 
functions on $X$, we have that  the function $F$ 
defined on $X\type k$ by
$$
 F(\ulx)= \prod_{\epsilon\in\{0,1\}^k}f_\epsilon(x_\epsilon)
$$
has the same integral under the measures $\mu\type k$ and $\nu_k$.

For every $x\in X$, the point $x\type k=(x,x,\dots,x)$ belongs to $X_k$. 
Since $(X_k,T\type k,T\type k_1,\dots, T\type k_k)$ 
is uniquely ergodic with invariant measure 
$\nu_k$, we have that
\begin{multline*}
\int F(\ulx)\,d\nu_k(\ulx)\\
 =\lim_{L\to+\infty}\frac 1L\sum_{\ell=0}^{L-1}\Bigl(
\lim_{M\to+\infty}\frac 1{M^{k-1}}\sum_{m_1,\dots 
m_{k-1}=0}^{M-1} \Bigl(
\lim_{N\to+\infty}\frac 1N\sum_{n=0}^{N-1}
\prod_{\epsilon\in\{0,1\}^k}
f_\epsilon(T^{n+\epsilon\cdot m+
\epsilon_k \ell}x)\Bigr)\Bigr)\ 
\end{multline*}
where $m=(m_1, \ldots, m_{k-1})$ and $\epsilon\cdot m = 
\epsilon_1m_1+\ldots+\epsilon_{k-1}m_{k-1}$.  
By unique ergodicity of $(X,T,\mu)$, this is equal to
$$
  \lim_{L\to+\infty}\frac 1L\sum_{\ell=0}^{L-1}\Bigl(
\lim_{M\to+\infty}\frac 1{M^{k-1}}\sum_{m_1,\dots 
m_{k-1}=0}^{M-1}
\int 
\prod_{\epsilon\in\{0,1\}^k}
f_\epsilon(T^{\epsilon\cdot m+\epsilon_k 
\ell}x)\,d\mu(x)\Bigr)\ .
$$
We write each $\epsilon\in\{0,1\}^k$ in the form 
$\eta 0$ or $\eta 1$ with $\eta\in\{0,1\}^{k-1}$, 
and this expression can be 
rewritten as 
$$
 \lim_{L\to+\infty}
\frac 1{L}\sum_{\ell=0}^{L-1}\Bigl(\lim_{M\to+\infty}
\frac 1{M^{k-1}} 
\sum_{m_1,\dots,m_{k-1}=0}^{M-1}
\int 
\prod_{\eta\in\{0,1\}^{k-1}}(f_{\eta 0}.T^\ell f_{\eta 1})
(T^{\eta\cdot m}x)\, d\mu(x)\Bigr)\ .
$$
By  unique ergodicity of $X_{k-1}$ under the 
transformations $T\type{k-1}$ and $T\type{k-1}_i$, $1\leq i\leq k-1$, 
and proceeding as above, we have that
this expression is equal to
$$
 \lim_{L\to+\infty}
\frac 1{L}\sum_{\ell=0}^{M-1}
\int \prod_{\eta\in\{0,1\}^{k-1}}(f_{\eta 0}.T^\ell f_{\eta 1})(x_\eta)\, 
d\nu_{k-1}(\ulx)\ .
$$
By the induction hypothesis, the integral remains unchanged 
when the measure $\mu\type{k-1}$ is substituted for $\nu_{k-1}$.
We rewrite this expression as
\begin{equation}
\label{eq:FFT}
 \lim_{L\to+\infty}
\frac 1{L}\sum_{\ell=0}^{L-1}\int F_0\,.\,  F_1\circ 
(T\type{k-1})^\ell\,d\mu\type{k-1} 
\end{equation}
where
$$
F_0(\ulx)=\prod_{\eta\in\{0,1\}^{k-1}}f_{\eta 0}(x_\eta)\text{ and }
F_1(\ulx)=\prod_{\eta\in\{0,1\}^{k-1}}f_{\eta 1}(x_\eta)\ .
$$
Let $\CI$ denotes the $T\type{k-1}$-invariant $\sigma$-algebra of the 
measure $\mu\type{k-1}$. The limit~\eqref{eq:FFT} is equal to
$$
 \int\E(F_0\mid\CI)\,\E(F_1\mid\CI)\,d\mu\type{k-1}\ .
$$
By the inductive definition of the measure $\mu\type k$ in~\cite{HK} (section 3), 
this is equal to
$$
 \int F_0(x_{\eta 0}\colon \eta\in\{0,1\}^{k-1})\,
F_1(x_{\eta 1}\colon \eta\in\{0,1\}^{k-1})\,
d\mu\type k(\ulx)
$$
and the function in the integral is just the function $F$.
\end{proof}
Recall that the measure $\mu\type k$ satisfies the inequality~\eqref{eq:HKCSG} 
of Section~\ref{sec:ergseminorms}. This 
can probably be proved directly for the measure $\nu_k$ but  does 
not seem obvious.

\subsection{The fibers}
Recall that for every $x\in X$,
$W_{k,x}=\{ \ulx\in X_k\colon x_\bzero=x\}$.

\begin{enumerate}
\enumitruc
\item\label{it:ergoWx}
For every $x\in X$, $W_{k,x}$ is uniquely ergodic under the 
transformations $T\type k_i$, $1\leq i\leq k$.
\end{enumerate}
\trucenumi
We write $\rho_x$ for the invariant measure of $W_{k,x}$.
\begin{enumerate}
\enumitruc
\item\label{it:invarWx}
For every $x\in X$ and $h\in G$, $\rho_{h.x}$ is the image of 
$\rho_x$ under the translation by $h\type k$.
\end{enumerate}
\trucenumi
\begin{proof}
Let $g$ be a lift of $x$ in $G$ and $\tilde\tau=g\tau g\inv$.

For $1\leq i\leq k$, we have that  $\tilde\tau\type k_i= g\type k \tau\type 
k_i(g\type k)\inv$ and all these elements commute and belong to $H_k$. 
For $1\leq i\leq k$, let $\tilde T\type k_i$ be the 
translation by $\tilde\tau\type k_i$.

We first show that the nilsystem $(W_k,\tilde T\type k_1,\dots,\tilde T\type k_k)$ 
is uniquely ergodic.
For each $i$, $\tilde\tau\type k_i(\tau\type k_i)\inv$ belongs to 
$H_k\cap (G_2)\type k$ and thus to $(H_k)_2$ by~\eqref{it:K2}.
Therefore, $\tilde\tau\type k_i$ and $\tau\type k_i$ have the same 
projection on the compact abelian group $H_k/(H_k)_2$. 
By~\eqref{it:Wkerg}, the action induced by  $\tau\type k_i$, 
$1\leq i\leq k$ on this group is ergodic.
The criteria given by Theorem~\ref{th:multiergo}  and property~\eqref{it:Hkconnected} give 
the announced unique ergodicity.

By~\eqref{it:WxHk}, we have that $g\type k.W_k=W_{k,x}$.
The map $\uly\mapsto g\type k.\uly$ mapping 
$(W_k,\tilde T\type k_1,\dots,\tilde T\type k_k)$ to 
$(W_{k,x}, T\type k_1,\dots,T\type k_k)$
 is an isomorphism of topological systems and thus the second of these 
system is uniquely ergodic. This proves~\eqref{it:ergoWx}.

We write $\rho$ for the Haar measure of the nilmanifold 
$W_k=H_k/\Theta_k$. Then $\rho$ is the invariant measure of 
$W_k$ and the above proof 
shows that for every $g\in G$, the invariant measure of $W_{k,x}$ is 
the image of $\rho$ under translation by $G\type k$. This 
immediately implies~\eqref{it:invarWx}.
\end{proof}

In fact,  $W_{k,x}$ can be given 
the structure of a nilmanifold, quotient of the group $H_k$ by the 
discrete cocompact group $g\type k\Theta(g\type k)\inv$, and the 
transformations $T\type k_i$ are translations on this nilmanifold.

\subsection{The case that $G$ is a $(k-1)$-step nilpotent}
Henceforth we assume that $G$ is a $(k-1)$-step nilpotent group. 

We show:
\begin{enumerate}
\enumitruc
\item  
Let $X_{k*}$ be the image of $\ulx\mapsto\ulx_*$ of $X_k$ under the 
projection $\ulx\mapsto\ulx_*$ mapping $X\type k$ to $X^{2^k-1}$.
There exists a smooth map $\Phi\colon X_{k*}\to X_k$ such that 
\begin{equation}\label{eq:Phi}
 X_k=\bigl\{(\Phi(\ulx_*),\ulx_*)\colon \ulx\in X_{k*}\bigr\}\ .
\end{equation}
\end{enumerate}
\trucenumi

Different proofs are given for the existence and continuity of $\Phi$ 
in~\cite{HK} and~\cite{GT}. 
The smoothness of $\Phi$ can be easily deduced from 
these proofs, but this property is not stated in these papers. For 
completeness, we give a short complete proof.

\begin{proof}
First we remark that the projection $X_k\to X_{k*}$ is one to one.
Indeed, let $\ulx$ and $\uly$ be two points of 
$X_k$ with the same projections. We lift them to two elements $\ulg$ 
and $\ulh$ of $G\type k$. All the coordinates of $\ulh\ulg\inv$ 
belong to $\Gamma$ except the first one, and by~\eqref{it:coordGamma} 
this coordinate also belongs to $\Gamma G_k=\Gamma$.  Thus $\ulx=\uly$.

Therefore the projection $X_k\to X_{k*}$ is a homeomorphism.
By composing the reciprocal of this  map with the projection 
$\ulx\mapsto x_\bzero$, we obtain a continuous map $\Phi\colon X_{k*}\to 
X$ satisfying~\eqref{eq:Phi}. We are left with showing that it is
smooth.

Let $G_*$ be the image of $G\type k$ in $G^{2^k-1}$ under the map 
$\ulg\mapsto\ulg_*$. By~\eqref{it:coordG}, the projection $G\type k\to 
G_*$ is one to one. 

We  check that $G_*$ is a closed subgroup of $G^{2^k-1}$.
Let $(\ulg_{*n})$ be a sequence in $G^{2^k-1}$ converging to some 
$\ulg_*$. For each $n$, there exists $g_{\bzero,n}\in G$ with 
$\ulg_n=(g_{\bzero,n},\ulg_{*n})\in G\type k_{k-1}$ and there exists
$\ulgamma_n\in\Gamma\type k\cap G\type k_{k-1}$ at a bounded distance 
from $\ulg_n$. All the  coordinates of $\ulgamma_n$, except 
$\gamma_\bzero$, are for all $n$ at a bounded distance from the unit.  
By passing to subsequences,  we can assume that they do not depend on 
$n$. By~\eqref{it:coordGamma}, $\ulgamma_n$ does not depend on $n$.
Therefore, $\ulg_n$ remains at a bounded distance from the unit and 
taking a subsequence we can assume that it converges to some $\ulg$, 
which belongs to $G\type k_{k-1}$ by~\eqref{it:Gkclosed}. Then the 
projection of $\ulg$ on $G_*$ is equal to $\ulg_*$.  Thus $\ulg$ 
belongs to $G_*$.

Now, the projection $G\type k_{k-1}\to G_*$ is a smooth bijective 
homomorphism between Lie groups.  Therefore it is a diffeomorphism. 
Since the 
projection $G\type k_{k-1}\to X_k$ has discrete kernel, it follows 
that the projection $X_k\to X_{k*}$ is a diffeomorphism and thus 
that $\Phi$ is smooth.
\end{proof}

We deduce: 
\begin{enumerate}
\enumitruc
\item
$\nnorm\cdot_k$ is a norm on $\CC(X)$.
\end{enumerate}

\begin{proof}
It suffices to show that if $f\in \CC(X)$ satisfies $\nnorm 
f_k=0$, then $f=0$.
By Proposition~\ref{prop:HKCSG},
if $f_\epsilon$, $\epsilon\in\{0,1\}^k_*$, 
are $2^k-1$ continuous functions on $X$, then
$$
 \int f(x_\bzero)\prod_{\epsilon\in 
\{0,1\}^k-*}f_\epsilon(x_\epsilon)\,d\mu\type k(\ulx)=0\ .
$$
By density, $\int f(x_0)F(\ulx_*)\,d\mu(\ulx)=0$ for every 
continuous function $F$ on $X_{k*}$. Taking $F=\bar f\circ\Phi$ where 
$\Phi$ is as in statement~\ref{it:Phi} of Theorem~\ref{th:summary},
property~\eqref{eq:Phi} of this 
function gives
$$
 0=\int f(x_0) \bar f(\Phi(x_*))\,d\mu\type k(\ulx)=\int |f(x_0)|^2\,d\mu\type k(\ulx)
=\int|f(x)|^2\,d\mu(x)
$$
because the projection of $\mu\type k$ on $X$ is $\mu$.
\end{proof}


\begin{thebibliography}{99999}

\bibitem[A]{auslander}
J.~Auslander. 
{\em Minimal flows and their extensions.}
North Holland Publishing Co, Amsterdam, 1988.


\bibitem[AGH]{AGH}
L.~Auslander, L.~Green and F.~Hahn.  Flows on homogeneous spaces. 
{\em Ann. Math. Studies} {\bf 53}, Princeton University Press, 1963. 


\bibitem[BFKO]{BFKO}
J.~Bourgain, H.~Furstenberg, Y.~Katznelson, D.~Ornstein.  
Appendix on return-time sequences.  {\em Inst. Hautes \'Etudes Sci.
Publ. Math.} {\bf 69} (1989), 42-45.

\bibitem[BFW]{BFW}
V.~Bergelson, H.~Furstenberg and B.~Weiss.
Piecewise-Bohr sets of integers and combinatorial number theory.  
Algorithms Combin. {\bf 26}, Springer, Berlin (2006), 
13-37. 

\bibitem[BHK]{BHK}
V.~Bergelson, B.~Host and B.~Kra, with an Appendix
by I.~Ruzsa.  Multiple recurrence
and nilsequences.  {\em Inventiones Math.} {\bf 160} (2005), 261-303.

\bibitem[BL]{BL}
A.~Leibman and V.~Bergelson. 
Distribution of values of bounded generalized polynomials.  
{\em Acta Math.} {\bf 198}  (2007), 155-230.


\bibitem[E]{E} R.~Ellis. 
{\em Lectures on topological dynamics.}  
W. A. Benjamin Inc., New York, 1969.

\bibitem[F]{F}
H.~Furstenberg.  
Ergodic behavior of diagonal measures and a theorem of Szemer\'edi
  on arithmetic progressions.
{\em J. d'Analyse Math.} {\bf 31}  (1977), 204-256.  

\bibitem[G]{Gowers} W.~T.~Gowers.  A new proof of Szemer\'edi's Theorem.
{\em Geom. Funct. Anal.} {\bf 11} (2001), 465-588.

\bibitem[GT1]{GT0} B.~Green and T.~Tao.  
The primes contain arbitrarily long arithmetic progressions. 
To appear, {\em Annals of Math}.
Available at: http://arxiv.org/abs/math/0404188

\bibitem[GT2]{GT} B.~Green and T.~Tao.  Linear equations in the primes.
To appear, {\em Annals of Math}.
Available at: http://arxiv.org/abs/math/0606088

\bibitem [GT3] {GT2}B.~Green and T.~Tao.  Quadratic uniformity of the 
M\"obius function.  To appear, {\em Annales de l'Institut Fourier}.
Available at: http://arxiv.org/abs/math/0606087

\bibitem[GT4]{GT3}  B.~Green and T.~Tao.
An inverse theorem for the Gowers $U^3$-norm, with applications.
To appear, {\em Proc. Edinburgh Math. Soc.}
Available at: http://arxiv.org/abs/math/0503014

\bibitem[HK1]{HK} B.~Host and B.~Kra.
Nonconventional ergodic averages and nilmanifolds.  
{\em Ann. Math.} {\bf 161} (2005), 397-488.

\bibitem [HK2]{HK2} B.~Host and B.~Kra.  Analysis of two step nilsequences. 
Submitted.  Available at: http://arxiv.org/abs/0709.3241

\bibitem[KN]{KN}
L.~Kuipers and H.~Niederreiter.
{\em Uniform distribution of sequences}.
John Wiley and Sons, New York, 1974.

\bibitem[Lei]{leibman}
A.~Leibman.  Pointwise convergence of ergodic averages for polynomial 
sequences of rotations of a nilmanifold.  {\em Erg. Th. \& Dynam. Sys.}
{\bf 25} (2005), 201-213.

\bibitem[Les1]{lesigne}
E.~Lesigne.  Sur une nil-vari\'et\'e, les parties minimales associ\'ees 
\`a une translation sont uniquement ergodiques.  {\em Erg. Th. \& Dynam. 
Sys.} {\bf 11} (1991), 379-391.

\bibitem[Les2]{lesigne2}
E.~Lesigne.  Spectre quasi-discret et th\'eor\`eme ergodique de 
Wiener-Wintner pour les polyn\^omes.  
{\em Erg. Th. \& Dynam. 
Sys.} {\bf 13} (1993), 767-784. 

\bibitem[Q]{Q}
M. Queffelec. \emph{Substitution Dynamical Systems --
Spectral Analysis}. Lecture Notes in Math.
{\bf 1294} Springer-Verlag, New York (1987).


\bibitem[WW]{WW}
N.~Wiener and A.~Wintner.  Harmonic analysis and ergodic theory. 
{\em Amer. J. Math.} {\bf 63} (1941), 415-426.
\end{thebibliography}
\end{document}